\documentclass[11pt,reqno]{amsart}

\usepackage[a4paper]{geometry}
\usepackage[headings]{fullpage}
\usepackage[onehalfspacing]{setspace}

\usepackage{graphicx}

\usepackage{mathrsfs}
\usepackage{amsfonts, amssymb, amsmath, amsthm, esint}
\usepackage{cite}
\usepackage{latexsym}
\usepackage{enumerate}
\usepackage{array}
\usepackage{color}
\usepackage{subcaption}
\usepackage{caption}
\usepackage{mdframed}

\usepackage{tcolorbox}

\usepackage{cleveref}

\usepackage{tikz}

\allowdisplaybreaks

\usepackage{enumitem}
\setlist[itemize]{leftmargin=2em}

\theoremstyle{definition}

\newtheorem{assump}{Assumption}

\theoremstyle{theorem}
\newtheorem{theorem}{Theorem}
\newtheorem{corollary}{Corollary}
\newtheorem{lemma}{Lemma}

\theoremstyle{remark}

\newenvironment{eqn*}[0]
{\begin{equation*}}{\end{equation*}\ignorespacesafterend}

\newenvironment{dis}[0]
{\begin{equation*}}{\end{equation*}\ignorespacesafterend}

\newenvironment{eqn}[0]
{\begin{equation}}{\end{equation}\ignorespacesafterend}

\newcommand{\mb}[1]{\mathbb{#1}}
\newcommand{\mc}[1]{\mathcal{#1}}

\newcommand{\tr}[1]{\textrm{#1}}

\newcommand{\what}[1]{\widehat{#1}}
\newcommand{\wtilde}[1]{\widetilde{#1}}

\newcommand{\pd}[0]{\partial}
\newcommand{\ol}[1]{\overline{#1}}

\newcommand{\emp}[0]{\emptyset}
\newcommand{\sus}[0]{\subset}

\newcommand{\vvep}[0]{\boldsymbol{\varepsilon}}
\newcommand{\vdel}[0]{\boldsymbol{\delta}}

\newcommand{\vsig}[0]{\boldsymbol{\sigma}}

\newcommand{\vpi}[0]{\boldsymbol{\pi}}

\newcommand{\vf}[0]{\boldsymbol{f}}

\newcommand{\vn}[0]{\boldsymbol{n}}

\newcommand{\vp}[0]{\boldsymbol{p}}
\newcommand{\vq}[0]{\boldsymbol{q}}
\newcommand{\vr}[0]{\boldsymbol{r}}
\newcommand{\vs}[0]{\boldsymbol{s}}
\newcommand{\vt}[0]{\boldsymbol{t}}
\newcommand{\vu}[0]{\boldsymbol{u}}
\newcommand{\vv}[0]{\boldsymbol{v}}
\newcommand{\vw}[0]{\boldsymbol{w}}
\newcommand{\vx}[0]{\boldsymbol{x}}

\newcommand{\zz}[0]{\boldsymbol{0}}

\newcommand{\vP}[0]{\boldsymbol{P}}

\newcommand{\vV}[0]{\boldsymbol{V}}
\newcommand{\vW}[0]{\boldsymbol{W}}
\newcommand{\vX}[0]{\boldsymbol{X}}

\newcommand*\diff{\mathop{}\!\mathrm{d}}

\newcommand{\inn}[1]{\left\langle #1\right\rangle}

\newcommand{\R}[0]{\mathbb{R}}

\newcommand{\N}[0]{\mathbb{N}}

\newcommand{\enorm}[1]{{\left\vert\kern-0.25ex\left\vert\kern-0.25ex\left\vert #1 
		\right\vert\kern-0.25ex\right\vert\kern-0.25ex\right\vert}}

\DeclareMathOperator{\rot}{rot}

\let\div\relax
\DeclareMathOperator{\div}{div}
\DeclareMathOperator{\vdiv}{\mathbf{div}}

\title[Lowest-Order VEMs for Linear Elasticity Problems]{Lowest-Order Virtual Element Methods for Linear Elasticity Problems}

\author{Do Y. Kwak \and Hyeokjoo Park}

\date{\today}


\thanks{Department of Mathematical Sciences, Korea Advanced Institute of Science and Technology, Daejeon, 305-701, Korea (kdy@kaist.ac.kr, hjpark235@kaist.ac.kr), This work is partially supported by NRF, contract No. 2021R1A2C1003340.}


\begin{document}

\begin{abstract}
We present two kinds of lowest-order virtual element methods for planar linear elasticity problems. For the first one we use the nonconforming virtual element method with a stabilizing term. It can be interpreted as a modification of the nonconforming Crouzeix-Raviart finite element method as suggested in \cite{MR1972650} to the virtual element method. For the second one we use the conforming virtual element for one component of the displacement vector and the nonconforming virtual element for the other. This approach can be seen as an extension of the idea of Kouhia and Stenberg suggested in \cite{MR1343077} to the virtual element method. We show that our proposed methods satisfy the discrete Korn's inequality. We also prove that the methods are convergent uniformly for the nearly incompressible case and the convergence rates are optimal. 
\end{abstract}
 
\maketitle

\section{Introduction}

We consider the following planar linear elasticity problem in a convex polygonal domain $\Omega\sus\R^2$: Given external force $\vf$, find the displacement field $\vu$ such that
\begin{eqn}\label{eqn:ModelProb}
\left\{\begin{array}{rl}
-\vdiv\left(2\mu\vvep(\vu) + \lambda(\div\vu)\vdel\right) = \vf & \tr{in} \ \Omega, \\
\left(2\mu\vvep(\vu) + \lambda(\div\vu)\vdel\right)\vn = \zz & \tr{on} \ \Gamma_N \sus \pd\Omega, \\
\vu = \zz & \tr{on} \ \Gamma_D = \pd\Omega - \Gamma_N.
\end{array}\right.
\end{eqn}
Here $\vn$ denotes the exterior unit vector normal to $\pd\Omega$, $\mu$ and $\lambda$ are the Lam\'e constants, $\vdel$ is the $2\times 2$ identity matrix and
\begin{dis}
\vvep(\vu) = \frac{1}{2}\left(\nabla\vu + (\nabla\vu)^{\perp}\right), \quad \div\vu = \frac{\pd u_1}{\pd x_1} + \frac{\pd u_2}{\pd x_2}.
\end{dis}
It is known that $\mu$ has positive lower and upper bounds and $0 < \lambda < \infty$. When the parameter $\lambda$ approaches to infinity, the problem \eqref{eqn:ModelProb} describes the behavior of nearly incompressible materials.

A challenging issue with developing numerical methods for this problem is the so-called locking phenomena, which may appear in the case of $\lambda\to\infty$ (i.e., the case of nearly incompressible materials). For instance, the piecewise linear conforming finite element method (FEM) may not converge as $\lambda\to\infty$. In order to obtain optimal convergence rates uniform with respect to $\lambda$ using the conforming FEM, the order of the finite element must be larger than $3$ \cite{MR1094947}. On the other hand, it was shown in \cite{MR1094947} that the nonconforming FEM with order $\geq 2$ has optimal convergence rates uniform with respect to $\lambda$. However, one cannot use the linear nonconforming element since the discrete Korn's inequality fails for this space, which means that the discrete bilinear form may not be coercive in general. 

Some researchers have developed some low-order FEMs while avoiding the difficulties mentioned above. For example, Hansbo and Larson \cite{MR1972650} used the linear nonconforming finite element with a stabilizing term to enforce the coerciveness. On the other hand, Kouhia and Stenberg \cite{MR1343077} proposed an element consisting of the linear conforming finite element for one component and the linear nonconforming finite element for the other. 

Meanwhile, the virtual element method (VEM) was recently introduced in \cite{beirao2013basic} as a generalization of the finite element method (FEM) to general polygonal or polyhedral meshes. In the VEM, the local discrete space on each mesh element consists of polynomials up to a given degree and some additional non-polynomial functions. In order to discretize continuous problems, the VEM only requires the knowledge of the degrees of freedom of the shape functions, such as values at mesh vertices, the moments on mesh edges/faces, or the moments on mesh polygons/polyhedrons, instead of knowing the shape functions explicitly. Moreover, the discrete space can be extended to high order in a straightforward way. Due to such advantages, The VEM has been successfully applied to various problems. For example, VEMs for general second-order elliptic problems were presented in \cite{MR3671497}. Some VEMs for the Stokes problem were introduced in \cite{MR3626409,MR4032861,MR4279152,MR3576570}. In \cite{MR3764431,MR3860122}, some VEMs for the magnetostatic problems were developed. For more thorough survey, we refer to \cite{MR3489088,MR3073346,MR3200242,MR3264352,MR3194807,MR3709049} and references therein. 

The VEM was also successfully applied to the linear elasticity problem. In \cite{MR3033033}, conforming virtual elements of order $\geq 2$ for the problem \eqref{eqn:ModelProb} were developed and it was shown that the convergence rates are optimal and uniform with respect to $\lambda$. In \cite{MR3881593}, the nonconforming VEMs with order $\geq 2$ were developed and it was shown that the convergence rates are optimal and uniform with respect to $\lambda$, but the lowest-order nonconforming VEM was developed only for the pure displacement problem, since the discrete Korn's inequality may fail for the lowest-order case. 

In this paper, we develop two kinds of lowest-order VEMs for the linear elasticity problem \eqref{eqn:ModelProb}. For the first one we use the lowest-order nonconforming virtual element with a stabilizing term. It can be interpreted as an extension of the method proposed by Hansbo and Larson \cite{MR1972650} to the virtual element method. For the second one we use the conforming virtual element for one component of the displacement field and the nonconforming virtual element for the other. This is similar to the element suggested by Kouhia and Stenberg \cite{MR1343077}. We then show that the proposed elements satisfy the discrete Korn's inequality. We also prove that the methods are convergent uniformly for the nearly incompressible case and the convergence rates are optimal under the regularity assumption of the solution. Here we only consider the case $\Gamma_N = \emp$ for convenience, but our proposed methods and their analysis can be easily extended to the general case $\Gamma_N\neq\emp$. 

The rest of our paper is organized as follows. In Section 2, we present model problem in a weak form, notations including mesh regularity, etc. In Section 3, we present the lowest-order nonconforming VEM with stabilizing term and prove optimal convergence. In Section 4, we present the Kouhia-Stenberg type VEM, and prove its convergence. In Section 5, we offer some numerical experiments to verify the performance of the proposed methods. Finally, conclusions are given in Section 6.

\section{Preliminaries}

Throughout this paper, we will use the usual Sobolev spaces $H^s(D)$, where $s \geq 0$ is an integer and $D$ is a bounded domain in $\R$ or $\R^2$. By convention, we note $H^0(D) = L^2(D)$. We denote by $\|\cdot\|_{s,D}$ and $|\cdot|_{s,D}$ the usual Sobolev norm and seminorm on $H^s(D)$, $(H^s(D))^2$, or $(H^s(D))^{2\times 2}$, respectively. We also denote $(\cdot,\cdot)_{0,D}$ the usual $L^2$-inner product on $L^2(D)$, $(L^2(D))^2$, or $(L^2(D))^{2\times 2}$. We also define
\begin{dis}
L_0^2(D) := \left\{q\in L^2(D) : \int_Dq\diff\vx = 0\right\}.
\end{dis}
For $s\geq 0$, we denote by $\mb{P}_s$ the space of polynomials of degree $\leq s$. 

\subsection{Model problem}

The linear elasticity problem \eqref{eqn:ModelProb} with $\Gamma_N = \emp$ has the following weak formulation: Given $\vf\in (L^2(\Omega))^2$, find $\vu\in (H_0^1(\Omega))^2$ such that
\begin{eqn}\label{eqn:ModelProbVar}
a(\vu,\vv) = (\vf,\vv)_{0,\Omega} \quad \forall \vv\in(H_0^1(\Omega))^2,
\end{eqn}
where
\begin{dis}
a(\vu,\vv) := 2\mu\int_\Omega\vvep(\vu):\vvep(\vv)\diff\vx + \lambda\int_{\Omega}\div\vu\div\vv\diff\vx.
\end{dis}
Note that the bilinear form $a(\cdot,\cdot)$ is bounded: there exists a positive constant $C$ independent of $\lambda$ such that
\begin{dis}
|a(\vu,\vv)| \leq C(1 + \lambda)|\vu|_{1,\Omega}|\vv|_{1,\Omega} \quad \forall \vu,\vv\in (H_0^1(\Omega))^2.
\end{dis}
Due to Korn's inequality \cite{MR2047078}, we obtain the ellipticity of $a(\cdot,\cdot)$: there exists a positive constant $C$ independent of $\lambda$ such that
\begin{dis}
C|\vv|_{1,\Omega}^2 \leq a(\vv,\vv) \quad \forall \vv\in (H_0^1(\Omega))^2.
\end{dis}
The boundedness and ellipticity of $a(\cdot,\cdot)$ shows that the problem \eqref{eqn:ModelProbVar} has a unique solution. Moreover, the solution $\vu$ of \eqref{eqn:ModelProbVar} satisfies the following regularity estimate \cite{MR1140646}: there exists a positive constant $C_{\Omega}$ depending only on $\Omega$ such that
\begin{eqn}\label{eqn:regularity}
\|\vu\|_{2,\Omega} + \lambda\|\div\vu\|_{1,\Omega} \leq C_{\Omega}\|\vf\|_{0,\Omega}.
\end{eqn}

\subsection{Mesh regularity}\label{subsec:MeshReg}

Let $\{\mc{P}_h\}_h$ be a sequence of decompositions of $\Omega$ into polygonal elements $K$ with maximum diameter $h$. Let $\mc{E}_h^i$ and $\mc{E}_h^{b}$ denote the set of all interior and boundary edges in $\mc{P}_h$, respectively. Similarly, let $\mc{V}_h^i$ and $\mc{V}_h^b$ be the set of all interior and boundary vertices in $\mc{P}_h$, respectively. We set $\mc{E}_h = \mc{E}_h^i\cup\mc{E}_h^b$ and $\mc{V}_h = \mc{V}_h^i\cup\mc{V}_h^b$. 

We assume that $\{\mc{P}_h\}_h$ satisfies the following regularity assumptions \cite{beirao2013basic,MR3507277,MR3709049}.
\begin{assump}\label{assum:MeshReg}
There exists $\rho > 0$ independent of $h$ such that
\begin{enumerate}[label=(\roman*)]
\item the decomposition $\mc{P}_h$ consists of a finite number of nonoverlapping polygonal elements;
\item for any $K\in\mc{P}_h$, the diameter of any edge of $K$ is larger than $\rho h_K$, where $h_K$ denotes the diameter $K$;
\item every element $K$ of $\mc{P}_h$ is star-shaped with respect to a ball with center $\vx_K$ and radius $\rho h_K$;
\item each element $K\in\mc{P}_h$ contains at least one interior vertex in $\mc{P}_h$.
\end{enumerate}
\end{assump}

Note that these assumptions imply the following properties \cite{MR3709049}:
\begin{itemize}
\item Every element $K\in\mc{P}_h$ has at most $N$ edges and vertices, where $N$ is independent of $h$.
\item For each element $K\in\mc{P}_h$, there is a triangular decomposition $\mc{T}^K$ obtained by connecting the vertices of $K$ to $\vx_K$ (see, for example, \Cref{fig:PolygonTriDecomp}), and the minimum angle of the triangular decomposition $\mc{T}^K$ is controlled by $\rho$.
\end{itemize}
For each $h$, we let
\begin{dis}
\mc{T}_h = \bigcup_{K\in\mc{P}_h}\mc{T}^K.
\end{dis}

\begin{figure}
\begin{center}
\begin{tikzpicture}[scale = 0.75]
\fill (0,0) circle [radius = 0.075];
\draw (1.5,0) -- (2,1) -- (0.5,2) -- (-1.5,1.25) -- (-2,0.25) -- (-1.5,-1.25) -- (-0.5,-1.5) -- (1,-1.25) -- (2,-0.75) -- (1.5,0);
\draw [dashed] (0,0) -- (1.5,0);
\draw [dashed] (0,0) -- (2,1);
\draw [dashed] (0,0) -- (0.5,2);
\draw [dashed] (0,0) -- (-1.5,1.25);
\draw [dashed] (0,0) -- (-2,0.25);
\draw [dashed] (0,0) -- (-1.5,-1.25);
\draw [dashed] (0,0) -- (-0.5,-1.5);
\draw [dashed] (0,0) -- (1,-1.25);
\draw [dashed] (0,0) -- (2,-0.75);
\draw [dashed] (0,0) circle [radius = 1];
\node at (-0.17,0.42) {$\vx_{K}$};
\node [right] at (1.5,-1.25) {$K$};
\end{tikzpicture}
\caption{}
\label{fig:PolygonTriDecomp}
\end{center}
\end{figure}

For each $K\in\mc{P}_h$, let $\vn_K$ and $\vt_K$ denote its exterior unit normal vector and counterclockwise tangential vector, respectively. For $e\in\mc{E}_h^i$, we define respectively $\vn_e$ and $\vt_e$ by a unit normal and tangential vector of $e$ with orientation fixed once and for all. For $e\in\mc{E}_h^b$, we define respectively $\vn_e$ and $\vt_e$ by a unit normal and tangential vector on $e$ in the outward and counterclockwise direction with respect to $\Omega$.

Let $e\in\mc{E}_h^i$ and let $K^-$ and $K^+$ be the polygons in $\mc{P}_h$ having $e$ as a common edge. For $v:\Omega\to\R$ satisfying $v|_{K^+}\in H^1(K^+)$ and $v|_{K^-}\in H^1(K^-)$, we define the jump of $v$ on $e$ by
\begin{dis}
[v]_e = v|_{K^+}(\vn_e\cdot\vn_{K^+}) + v|_{K^-}(\vn_e\cdot\vn_{K^-}).
\end{dis}
If $e\in\mc{E}_h^b$, we define $[v]_e = v|_e$. Analogously, we define $[\vv]_e = ([v_1]_e,[v_2]_e)$ for $\vv:\Omega\to\R^2$ with $\vv = (v_1,v_2)$ satisfying $\vv|_{K^+}\in (H^1(K^+))^2$ and $\vv|_{K^-}\in(H^1(K^-))^2$. 

We let $C$ denote a generic positive constant independent of the Lam\'e constant $\lambda$ and the mesh parameter $h$, not necessarily the same in each occurrence. 

Given a decomposition $\mc{P}$ of $\Omega$ into a finite number of non-overlapping polygonal elements, we define the broken Sobolev space
\begin{dis}
H^1(\Omega;\mc{P}) = \left\{v\in L^2(\Omega) : v|_K\in H^1(K) \ \forall K\in\mc{P}\right\}.
\end{dis}
We also define the broken $H^1$-seminorm on the space $H^1(\Omega;\mc{P})$ or $(H^1(\Omega;\mc{P}))^2$ as follows:
\begin{dis}
|\cdot|_{H^1(\Omega;\mc{P})}^2 = \sum_{K\in\mc{P}}|\cdot|_{1,K}^2. 
\end{dis}
In particular, if $\mc{P} = \mc{P}_h$ then we simply write $|\cdot|_{1,h} = |\cdot|_{H^1(\Omega;\mc{P}_h)}$.

For $\vv\in (H^1(\Omega;\mc{T}_h))^2$, we define $\nabla_h\vv$ by $(\nabla_h\vv)|_T = \nabla(\vv|_T)$ for each $T\in\mc{T}_h$. Analogous definitions hold for $\vvep_h$, $\rot_h$, and $\div_h$. Here the operator $\rot$ is defined by $\rot\vv = \frac{\pd v_2}{\pd x_1} - \frac{\pd v_1}{\pd x_2}$ for a field $\vv = (v_1,v_2)$. 

For convenience, we define the local bilinear form $a^K:(H^1(K))^2\times(H^1(K))^2\to\R$ on each element $K\in\mc{P}_h$ by $a^K = a_{\mu}^K + a_{\lambda}^K$ where
\begin{dis}
a_{\mu}^K(\vu,\vv) = 2\mu\int_K\vvep(\vu):\vvep(\vv)\diff\vx, \quad a_{\lambda}^K(\vu,\vv) = \lambda\int_K\div\vu\div\vv\diff\vx.
\end{dis}


\section{Lowest-Order Nonconforming VEM with Stabilizing Term}\label{sec:NCVEMSTAB}

In this section, we present the lowest-order nonconforming VEM for the problem \eqref{eqn:ModelProbVar}. 

\subsection{Lowest-order nonconforming virtual element space}

Let $K$ be a polygon satisfying the regularity assumptions (ii) and (iii) in \Cref{assum:MeshReg}. We first define an auxiliary local space $\wtilde{\vV}_h(K) = (\wtilde{V}_h(K))^2$, where 
\begin{dis}
\wtilde{V}_h(K) = \left\{v\in H^1(K) : \Delta v = 0, \ (\vn_K\cdot\nabla v)|_e\in \mb{P}_1(e) \ \forall e\sus \pd K\right\}.
\end{dis}
We also define a projection operator $\Pi_h^K:\wtilde{\vV}_h(K)\to (\mb{P}_1(K))^2$ as the solution of
\begin{dis}
\left\{\begin{array}{l}
\int_K\vvep(\Pi_h^K\vv):\vvep(\vq)\diff\vx = \int_K\vvep(\vv):\vvep(\vq)\diff\vx \quad \forall \vq\in (\mb{P}_1(K))^2, \\
\int_K\rot\Pi_h^K\vv\diff\vx = \int_K\rot\vv\diff\vx, \\
\int_{\pd K}\Pi_h^K\vv\diff s = \int_{\pd K}\vv\diff s,
\end{array}\right.
\end{dis}
for $\vv\in\wtilde{\vV}_h(K)$. Note that
\begin{eqnarray*}
\int_K\vvep(\vv):\vvep(\vq)\diff\vx & = & \int_{\pd K}(\vvep(\vq)\vn_K)\cdot\vv\diff s \quad \forall \vq\in (\mb{P}_1(K))^2, \\
\int_K\rot\vv\diff\vx & = & \int_{\pd K}\vv\cdot\vt_K\diff s,
\end{eqnarray*}
for any $\vv\in\wtilde{\vV}_h(K)$. Thus $\Pi_h^K\vv$ is computable from the degrees of freedom
\begin{eqn}\label{eqn:NCVEMDOF}
\tr{the moments} \ \frac{1}{|e|}\int_e\vv\diff s, \quad \forall e\sus \pd K.
\end{eqn}
Moreover, $\Pi_h^K\vq = \vq$ for any $\vq\in(\mb{P}_1(K))^2$. We then define the local virtual element space
\begin{eqn}\label{eqn:NCVEMLocal}
\vV_h(K) = \left\{\vv\in\wtilde{\vV}_h(K) : \int_e\vv\cdot\vq\diff s = \int_e\Pi_h^K\vv\cdot\vq\diff s \ \forall \vq\in (\mb{P}_1^*(e))^2 \ \forall e\sus\pd K\right\},
\end{eqn}
where $\mb{P}_1^*(e)$ denotes the subspace of $\mb{P}_1(e)$ that is $L^2(e)$-orthogonal to $\mb{P}_0(e)$. It is not difficult to show that $(\mb{P}_1(K))^2\sus \vV_h(K)$ and the (local) degrees of freedom \eqref{eqn:NCVEMDOF} are unisolvent for $\vV_h(K)$. 

The global virtual element space $\vV_h$ is defined by
\begin{dis}
\vV_h = \left\{\vv\in (L^2(\Omega))^2 : \vv|_K\in\vV_h(K) \ \forall K\in\mc{P}_h, \ \int_e[\vv]_e\diff s = \zz \ \forall e\in\mc{E}_h\right\}.
\end{dis}
It is easy to see that the following degrees of freedom are unisolvent for $\vV_h$:
\begin{dis}
\tr{the moments} \ \frac{1}{|e|}\int_e\vv\diff s, \quad \forall e\in\mc{E}_h^i.
\end{dis}

Given $\vv\in (H_0^1(\Omega))^2$, we denote by $I_h\vv$ the global interpolant of $\vv$, that is, $I_h\vv$ is defined by the unique function in $\vV_h$ such that $\chi_i(\vv - I_h\vv) = 0$ for any $i = 1,2,\cdots,\dim\vV_h$, where $\chi_i$ is the operator that associates the $i$-th degree of freedom of $\vV_h$. Then the following lemma holds (see (3.16) in \cite{MR3507277}).

\begin{lemma}\label{lem:NCInterpError}
Let $I_h:(H_0^1(\Omega))^2\to \vV_h$ be the interpolation operator as defined above. There exists a positive constant $C$ independent of $h$ such that for any $\vv\in (H_0^1(\Omega)\cap H^2(\Omega))^2$ and any $K\in\mc{P}_h$,
\begin{dis}
\|\vv - I_h\vv\|_{0,K} + h_K|\vv - I_h\vv|_{1,K} \leq Ch_K^2|\vv|_{2,K}.
\end{dis}
\end{lemma}

\subsection{Discrete problem}

Let $\Pi_0^K:L^2(K)\to \mb{P}_0(K)$ be the $L^2$-projection operator. Let $S^K:\vV_h(K)\times\vV_h(K)\to\R$ be a bilinear form such that
\begin{dis}
S^K(\vu,\vv) = \sum_{i=1}^{\dim\vV_h(K)}\chi_i(\vu)\chi_i(\vv),
\end{dis}
where $\chi_i$ is the operator associated with the $i$-th local degrees of freedom. We then define the local bilinear forms $a_{\mu,h}^K$, $a_{\lambda,h}^K$, and $a_h^K$ on $\vV_h(K)$ by
\begin{eqnarray*}
a_h^K(\vu,\vv) & = & a_{\mu,h}^K(\vu,\vv) + a_{\lambda,h}^K(\vu,\vv), \quad \tr{where} \\
a_{\mu,h}^K(\vu,\vv) & = & 2\mu\int_K\vvep(\Pi_h^K\vu):\vvep(\Pi_h^K\vv)\diff\vx+ S^K(\vu - \Pi_h^K\vu, \vv - \Pi_h^K\vv), \\
a_{\lambda,h}^K(\vu,\vv) & = & \lambda\int_K(\Pi_0^{K}\div\vu)(\Pi_0^{K}\div\vv)\diff\vx.
\end{eqnarray*}
for $\vu,\vv\in\vV_h(K)$. Following the arguments in \cite{beirao2013basic}, it is easy to show that the bilinear form $a_h^K(\cdot,\cdot)$ satisfies the consistency and the stability:
\begin{itemize}
\item (Consistency) $a_h^K(\vp,\vv) = a^K(\vp,\vv)$ for any $\vv\in\vV_h(K)$ and $\vp\in(\mb{P}_1(K))^2$;
\item (Stability) there exist two positive constant $c_*$ and $c^*$, independent of $h$ and of $K$, such that $c_*a_{\mu}^K(\vv,\vv) \leq a_{\mu,h}^K(\vv,\vv) \leq c^*a_{\mu}^K(\vv,\vv)$ for any $\vv\in\vV(K)$.
\end{itemize}

We next define the global discrete bilinear form $a_h:\vV_h\times\vV_h\to\R$. Note that, however, the bilinear form
\begin{dis}
\sum_{K\in\mc{P}_h}a_h^K(\vu_h,\vv_h), \quad \vu_h,\vv_h\in\vV_h,
\end{dis}
is not elliptic with respect to $|\cdot|_{1,h}$ since the functions in $\vV_h$ do not satisfy Korn's inequality in general. To avoid this, we will add a stabilizing term $J_h(\cdot,\cdot)$ as in \cite{MR1972650}. We define $J_h:\vV_h\times\vV_h\to\R$ by
\begin{dis}
J_h(\vu_h,\vv_h) = \frac{\gamma}{h}\sum_{e\in\mc{E}_h^i}\int_e\vpi_e[\vu_h]_e\cdot\vpi_e[\vv_h]_e\diff s,
\end{dis}
where $\vpi_e$ denotes the $L^2$-projection from $(L^2(e))^2$ onto $(\mb{P}_1(e))^2$ for each $e\in\mc{E}_h^i$ and $\gamma$ is a fixed positive constant. Note that, due to \eqref{eqn:NCVEMLocal}, if $e\in\mc{E}_h^i$ is a common edge of the elements $K^+,K^-\in\mc{P}_h$, then 
\begin{eqnarray*}
\int_e\vpi_e[\vv_h]_e\cdot\vq\diff s & = & \int_e((\vv_h|_{K^+})\cdot\vq)(\vn_{K^+}\cdot\vn_e)\diff s + \int_e((\vv_h|_{K^-})\cdot\vq)(\vn_{K^-}\cdot\vn_e)\diff s \\
& = & \int_e((\Pi_h^K(\vv_h|_{K^+}))\cdot\vq)(\vn_{K^+}\cdot\vn_e)\diff s + \int_e((\Pi_h^K(\vv_h|_{K^-}))\cdot\vq)(\vn_{K^-}\cdot\vn_e)\diff s
\end{eqnarray*}
for any $\vq\in(\mb{P}_1^*(e))^2$ and $\vv_h\in\vV_h$. Thus $J_h(\vu_h,\vv_h)$ is computable using only the degrees of freedom of $\vu_h,\vv_h\in\vV_h$.

Now the global discrete bilinear form $a_h:\vV_h\times\vV_h\to\R$ is defined by
\begin{dis}
a_h(\vu_h,\vv_h) = \sum_{K\in\mc{P}_h}a_h^K(\vu_h,\vv_h) + J_h(\vu_h,\vv_h), \quad \forall \vu_h,\vv_h\in\vV_h.
\end{dis}
We will later show that $a_h(\cdot,\cdot)$ is elliptic with respect to $|\cdot|_{1,h}$.

We next construct the discrete loading term. We define $\vf_h$ on each element $K\in\mc{P}_h$ as the $(L^2(K))^2$-projection of $\vf$ on the space of piecewise constant, that is,
\begin{dis}
\vf_h|_K = \frac{1}{|K|}\int_K\vf\diff\vx, \quad \forall K\in\mc{P}_h,
\end{dis}
where $|K|$ denotes the area of $K$. We then define the discrete loading term $\inn{\vf_h,\cdot}$ as follows:
\begin{dis}
\inn{\vf_h,\vv_h} = \sum_{K\in\mc{P}_h}\int_K\vf_h\cdot\what{\vv}_h\diff\vx, \quad \forall \vv_h\in\vV_h
\end{dis}
where $\hat{\vv}_h$ is defined by
\begin{dis}
\what{\vv}_h|_K = \frac{1}{N_{K}}\sum_{\substack{e\in\mc{E}_h\\ e\sus \pd K}}\frac{1}{|e|}\int_e\vv_h\diff s, \quad \forall K\in\mc{P}_h.
\end{dis}
Then the following lemma for the approximation of the loading term $(\vf,\cdot)_{0,\Omega}$ can be found in \cite{beirao2013basic,MR3507277}.

\begin{lemma}\label{lem:NCloading}
Suppose that $\vf\in (L^2(\Omega))^2$. Then there exists a positive constant $C$ independent of $h$ such that
\begin{dis}
\left|\inn{\vf_h,\vv_h} - (\vf,\vv_h)_{0,\Omega}\right| \leq Ch\|\vf\|_{0,\Omega}|\vv_h|_{1,h} \quad \forall \vv_h \in\vV_h.
\end{dis}
\end{lemma}

With the above preparations, we state the following virtual element discretization of the problem \eqref{eqn:ModelProbVar}: Find $\vu_h\in\vV_h$ such that
\begin{eqn}\label{eqn:DiscreteProbNCS}
a_h(\vu_h,\vv_h) = \inn{\vf_h,\vv_h} \quad \forall \vv_h\in\vV_h.
\end{eqn}

\subsection{Error analysis}

It is well-known that the following approximation property holds \cite{MR2373954}.

\begin{lemma}\label{lem:ApproxPoly}
Let $K\in\mc{P}_h$. For any $\vv\in (H^2(K))^2$, there exists $\vv_{\pi}\in(\mb{P}_1)^2$ such that
\begin{dis}
\|\vv - \vv_{\pi}\|_{0,K} + h_K|\vv - \vv_{\pi}|_{1,K} \leq Ch_K^2|\vv|_{2,K},
\end{dis}
where $C$ is a positive constant depending only on $\rho$.
\end{lemma}

We first check the existence and uniqueness of the solution of the discrete problem \eqref{eqn:DiscreteProbNCS}. According to the result in \cite{MR2047078}, there exists a positive constant $C$ independent of $h$ such that
\begin{dis}
|\vv_h|_{1,h}^2 \leq C\left(\|\vvep_h(\vv_h)\|_{0,\Omega}^2 + \left|\sum_{K\in\mc{P}_h}\int_K\rot\vv_h\diff\vx\right| + \sum_{e\in\mc{E}_h^i}\frac{1}{|e|}\|\vpi_e[\vv_h]_e\|_{0,e}^2\right).
\end{dis}
Since $\int_e[\vv_h]_e\diff s = \zz$ for any $\vv_h\in\vV_h$ and $e\in\mc{E}_h$, 
\begin{dis}
\sum_{K\in\mc{P}_h}\int_K\rot\vv_h\diff\vx = \sum_{K\in\mc{P}_h}\int_{\pd K}\vv_h\cdot\vt_K\diff s = \sum_{e\in\mc{E}_h}\int_e[\vv_h]_e\cdot\vt_e\diff s = 0.
\end{dis}
Therefore we deduce that there exists a positive constant $C$ independent of $h$ such that
\begin{eqn}\label{eqn:NCVEMKorn}
|\vv_h|_{1,h}^2 \leq C\left(\|\vvep_h(\vv_h)\|_{0,\Omega}^2 + \sum_{e\in\mc{E}_h^i}\frac{1}{|e|}\|\vpi_e[\vv_h]_e\|_{0,e}^2\right) \leq Ca_h(\vv_h,\vv_h) \quad \forall \vv_h\in\vV_h.
\end{eqn}
This inequality shows that the discrete bilinear form $a_h(\cdot,\cdot)$ is elliptic on $\vV_h$, and hence the discrete problem \eqref{eqn:DiscreteProbNCS} has a unique solution.

We next prove the following convergence theorem.

\begin{theorem}\label{thm:NCVEMconv}
Suppose that $\vf\in (L^2(\Omega))^2$ and $\vu\in (H^2(\Omega)\cap H_0^1(\Omega))^2$ is the solution of \eqref{eqn:ModelProbVar}. Let $\vu_h\in\vV_h$ be the unique solution of the discrete problem \eqref{eqn:DiscreteProbNCS}. Then
\begin{dis}
|\vu - \vu_h|_{1,h} \leq Ch\|\vf\|_{0,\Omega},
\end{dis}
where $C$ is a positive constant independent of $h$ and the Lam\'e constant $\lambda$. 
\end{theorem}

\begin{proof}
Let $\vu_{\pi}$ be the approximation in \Cref{lem:ApproxPoly}, $\vu_I = I_h\vu$, and $\vdel_h = \vu_h - \vu_I$. Define a norm $\enorm{\cdot}$ on $\vV_h$ by $\enorm{\vv_h}^2 := a_h(\vv_h,\vv_h)$. Using the consistency of $a_h(\cdot,\cdot)$,
\begin{align}
\enorm{\vdel_h}^2 = &\ a_h(\vdel_h,\vdel_h) = a_h(\vu_h,\vdel_h) - a_h(\vu_I,\vdel_h) = \inn{\vf_h,\vdel_h} - \sum_{K\in\mc{P}_h}a_h^K(\vu_I,\vdel_h) - J_h(\vu_I,\vdel_h) \nonumber\\
= &\ \inn{\vf_h,\vdel_h} - \sum_{K\in\mc{P}_h}\left(a_h^K(\vu_I - \vu_{\pi},\vdel_h) + a_h^K(\vu_{\pi},\vdel_h)\right) - J_h(\vu_I,\vdel_h) \nonumber\\
= &\ \inn{\vf_h,\vdel_h} - \sum_{K\in\mc{P}_h}\left(a_h^K(\vu_I - \vu_{\pi},\vdel_h) + a^K(\vu_{\pi},\vdel_h)\right) - J_h(\vu_I,\vdel_h) \nonumber\\
= &\ \inn{\vf_h,\vdel_h} - \sum_{K\in\mc{P}_h}\left(a_h^K(\vu_I - \vu_{\pi},\vdel_h) + a^K(\vu_{\pi} - \vu,\vdel_h)\right) - \sum_{K\in\mc{P}_h}a^K(\vu,\vdel_h) - J_h(\vu_I,\vdel_h) \nonumber\\
= &\ \inn{\vf_h,\vdel_h} - (\vf,\vdel_h)_{0,\Omega} - \sum_{K\in\mc{P}_h}\left(a_h^K(\vu_I - \vu_{\pi},\vdel_h) + a^K(\vu_{\pi} - \vu,\vdel_h)\right) \nonumber\\
& - \sum_{K\in\mc{P}_h}a^K(\vu,\vdel_h) + (\vf,\vdel_h)_{0,\Omega} - J_h(\vu_I,\vdel_h). \label{eqn:NCError00}
\end{align}
Let
\begin{align*}
& T_1 := \inn{\vf_h,\vdel_h} - (\vf,\vdel_h)_{0,\Omega}, \qquad T_2 := \sum_{K\in\mc{P}_h}\left(a_h^K(\vu_I - \vu_{\pi},\vdel_h) + a^K(\vu_{\pi} - \vu,\vdel_h)\right), \\
& T_3 := \sum_{K\in\mc{P}_h}a^K(\vu,\vdel_h) - (\vf,\vdel_h)_{0,\Omega}, \qquad T_4 := J_h(\vu_I,\vdel_h).
\end{align*}
Note that, by \eqref{eqn:NCVEMKorn},
\begin{eqn}\label{eqn:NCVEMEnorm}
C|\vv_h|_{1,h} \leq \enorm{\vv_h} \quad \forall \vv_h\in\vV_h.
\end{eqn}
By \Cref{lem:NCloading} and \eqref{eqn:NCVEMEnorm}, 
\begin{eqn}
\left|T_1\right| \leq Ch\|\vf\|_{0,\Omega}\enorm{\vdel_h}. \label{eqn:NCError01}
\end{eqn}
Note that
\begin{eqnarray*}
&& a_h^K(\vu_I - \vu_{\pi},\vdel_h) + a^K(\vu_{\pi} - \vu,\vdel_h) \\
& = & a_{\mu,h}^K(\vu_I - \vu_{\pi},\vdel_h) + a_{\mu}^K(\vu_{\pi} - \vu,\vdel_h) + a_{\lambda,h}^K(\vu_I - \vu_{\pi},\vdel_h) + a_{\lambda}^K(\vu_{\pi} - \vu,\vdel_h).
\end{eqnarray*}
Since $\vu_I$ has the same degrees of freedom with $\vu$,
\begin{dis}
\int_Kq\div\vu_I\diff\vx = \int_{\pd K}q\vu_I\cdot\vn_K\diff\vx = \int_{\pd K}q\vu\cdot\vn_K\diff\vx = \int_Kq\div\vu\diff\vx
\end{dis}
for any $q\in\mb{P}_0(K)$ and any $K\in\mc{P}_h$. Then we have
\begin{eqnarray*}
&& a_{\lambda,h}^K(\vu_I - \vu_{\pi},\vdel_h) + a_{\lambda}^K(\vu_{\pi} - \vu,\vdel_h) \\
& = & \lambda(\Pi_0^{K}\div \vu_I - \Pi_0^{K}\div\vu_{\pi},\Pi_0^{K}\div\vdel_h)_{0,K} + \lambda(\div\vu_{\pi} - \div\vu,\div\vdel_h)_{0,K} \\
& = & \lambda(\Pi_0^{K}\div \vu_I,\Pi_0^{K}\div\vdel_h)_{0,K} - \lambda(\div\vu,\div\vdel_h)_{0,K} \\
& = & \lambda(\div \vu_I,\Pi_0^{K}\div\vdel_h)_{0,K} - \lambda(\div\vu,\div\vdel_h)_{0,K} \\
& = & \lambda(\div \vu,\Pi_0^{K}\div\vdel_h)_{0,K} - \lambda(\div\vu,\div\vdel_h)_{0,K} \\
& = & \lambda(\Pi_0^{K}\div\vu,\div\vdel_h)_{0,K} - \lambda(\div\vu,\div\vdel_h)_{0,K} \\
& = & \lambda(\Pi_0^{K}\div\vu - \div\vu,\div\vdel_h)_{0,K}
\end{eqnarray*}
for any $K\in\mc{P}_h$. Using \Cref{lem:ApproxPoly}, \Cref{lem:NCInterpError} and \eqref{eqn:NCVEMEnorm} we obtain
\begin{eqnarray}
|T_2| & \leq & C\sum_{K\in\mc{P}_h}\left(|\vu_I - \vu_{\pi}|_{1,K} + |\vu - \vu_{\pi}|_{1,K} + \lambda|\Pi_0^{K}\div\vu - \div\vu|_{0,K}\right)|\vdel_h|_{1,K} \nonumber\\
& \leq & C\sum_{K\in\mc{P}_h}h_K\left(|\vu|_{2,K} + \lambda|\div\vu|_{1,K}\right)|\vdel_h|_{1,K} \nonumber\\
& \leq & Ch\left(|\vu|_{2,\Omega} + \lambda|\div\vu|_{1,\Omega}\right)\enorm{\vdel_h}. \label{eqn:NCError02}
\end{eqnarray}
Integrating by parts we obtain
\begin{dis}
T_3 = \sum_{K\in\mc{P}_h}\int_{\pd K}\left(\vsig(\vu)\vn_K\right)\cdot\vdel_h\diff s = \sum_{e\in\mc{E}_h}\int_e\vsig(\vu)\vn_e\cdot[\vdel_h]_e\diff s.
\end{dis}
For each $e\in\mc{E}_h$, let $\vP_e^0:(L^2(e))^2\to (\mb{P}_0(e))^2$ be the $L^2$-orthogonal projection operator. Then, since 
\begin{eqn}
\int_e[\vv_h]_e\diff s = 0 \quad \forall e\in\mc{E}_h, \ \forall \vv_h\in\vV_h, \label{eqn:NCDiscreteJump}
\end{eqn}
we obtain 
\begin{eqnarray*}
\int_e\vsig(\vu)\vn_e\cdot[\vdel_h]_e\diff s & = & \int_e\left(\vsig(\vu)\vn_e - \vP_e^0\vsig(\vu)\vn_e\right)\cdot[\vdel_h]_e\diff s \\
& = & \int_e\left(\vsig(\vu)\vn_e - \vP_e^0\vsig(\vu)\vn_e\right)\cdot[\vdel_h - \vP_e^0\vdel_h]_e\diff s
\end{eqnarray*}
From the classical arguement in \cite{MR343661}, if $e\in\mc{E}_h^i$ and $e$ is a common edge of two elements $K^+$ and $K^-$ in $\mc{P}_h$, then
\begin{eqnarray*}
\|\vsig(\vu)\vn_e - \vP_e^0\vsig(\vu)\vn_e\|_{0,e} & \leq & Ch_e^{1/2}\|\vsig(\vu)\|_{1,K^+\cup K^-}, \\
\|[\vdel_h - \vP_e^0\vdel_h]_e\|_{0,e} & \leq & Ch^{1/2}\left(|\vdel_h|_{1,K^+}^2 + |\vdel_h|_{1,K^-}^2\right)^{1/2}.
\end{eqnarray*}
If $e\in\mc{E}_h^b$, then \eqref{eqn:NCDiscreteJump} implies that $[\vdel_h - \vP_e^0\vdel_h]_e = 0$. Thus, using \eqref{eqn:NCVEMEnorm},
\begin{eqn}
\left|T_3\right| \leq Ch\left(\|\vu\|_{2,\Omega} + \lambda\|\div\vu\|_{1,\Omega}\right)\enorm{\vdel_h}\label{eqn:NCError03}
\end{eqn}
Since $\vu\in (H_0^1(\Omega)\cap H^2(\Omega))^2$, $[\vu]_e = 0$ for any $e\in\mc{E}_h$ and so
\begin{dis}
T_4 = J_h(\vu_I,\vdel_h) = J_h(\vu_I - \vu, \vdel_h).
\end{dis}
Then
\begin{eqnarray*}
|J_h(\vu_I - \vu, \vdel_h)| & \leq & \frac{\gamma}{h}\sum_{e\in\mc{E}_h^i}\|[\vu - \vu_I]_e\|_{0,e}\|[\vdel_h]_e\|_{0,e} \\
& \leq & \left(\frac{\gamma}{h}\sum_{e\in\mc{E}_h^i}\|[\vu - \vu_I]_e\|_{0,e}^2\right)^{1/2}\left(\frac{\gamma}{h}\sum_{e\in\mc{E}_h^i}\|[\vdel_h]_e\|_{0,e}\right)^{1/2} \\
& \leq & \left(\frac{\gamma}{h}\sum_{e\in\mc{E}_h^i}\|[\vu - \vu_I]_e\|_{0,e}^2\right)^{1/2}\enorm{\vdel_h}.
\end{eqnarray*}
Let $e\in\mc{E}_h^i$ and assume that $e$ is a common edge of two elements $K_1$ and $K_2$ in $\mc{P}_h$. From the trace theorem with scaling and \Cref{lem:NCInterpError}, we obtain
\begin{eqnarray*}
\|[\vu - \vu_I]_e\|_{0,e}^2 \leq C\sum_{i=1}^2\left(h|\vu - \vu_I|_{1,K_i}^2 + h^{-1}\|\vu - \vu_I\|_{0,K_i}^2\right) \leq Ch^3\left(|\vu|_{2,K_1}^2 + |\vu|_{2,K_2}^2\right).
\end{eqnarray*}
Thus
\begin{eqn}
|T_4| = |J_h(\vu_I - \vu, \vdel_h)| \leq C\left(\frac{\gamma}{h}\sum_{K\in\mc{P}_h}h^3|\vu|_{2,K}^2\right)^{1/2}\enorm{\vdel_h} \leq Ch|\vu|_{2,\Omega}\enorm{\vdel_h}. \label{eqn:NCError04}
\end{eqn}
Now combining the results \eqref{eqn:NCError00}, \eqref{eqn:NCError01}, \eqref{eqn:NCError02}, \eqref{eqn:NCError03}, and \eqref{eqn:NCError04}, we obtain
\begin{dis}
\enorm{\vdel_h}^2 \leq Ch\left(\|\vu\|_{2,\Omega} + \lambda\|\div\vu\|_{1,\Omega} + \|\vf\|_{0,\Omega}\|\right)\enorm{\vdel_h}.
\end{dis}
Finally, using the regularity estimate \eqref{eqn:regularity} and \eqref{eqn:NCVEMEnorm},
\begin{dis}
|\vu - \vu_h|_{1,h} \leq |\vu - \vu_I|_{1,h} + C\enorm{\vdel_h} \leq Ch|\vu|_{2,\Omega} + C\enorm{\vdel_h} \leq Ch\|\vf\|_{0,\Omega}.
\end{dis}
This concludes the proof of the theorem.
\end{proof}


\section{Kouhia-Stenberg type VEM}\label{sec:KSVEM}

In this section, we present the Kouhia-Stenberg type VEM for the problem \eqref{eqn:ModelProbVar}. 

\subsection{Kouhia-Stenberg type virtual element space}

Let $K$ be a polygon satisfying the regularity assumptions (ii) and (iii) in \Cref{assum:MeshReg}. We first introduce an auxiliary space
\begin{dis}
B(\pd K) = \left\{g\in C^0(\pd K) : g|_e\in \mb{P}_1(e) \ \forall e\sus \pd K\right\}.
\end{dis} 
Then the local conforming and nonconforming virtual element spaces are defined as follows \cite{beirao2013basic,MR3507277}:
\begin{eqnarray*}
V_{c}(K) & = & \left\{v\in H^1(K) : \Delta v = 0 \ \tr{in} \ K \ \tr{and} \ v|_{\pd K}\in B(\pd K) \right\}, \\
V_{nc}(K) & = & \left\{v\in H^1(K) : \Delta v = 0 \ \tr{in} \ K \ \tr{and} \ (\vn_K\cdot\nabla v)|_e \in \mb{P}_0(e) \ \forall e\sus\pd K\right\}.
\end{eqnarray*}
The conforming and nonconforming global virtual element spaces are defined by
\begin{eqnarray*}
V_{h,c} & = & \left\{v\in H_0^1(\Omega) : v|_K\in V_c(K) \ \forall K\in\mc{P}_h\right\}, \\
V_{h,nc} & = & \left\{v\in L^2(\Omega) : v|_K\in V_{nc}(K) \ \forall K\in\mc{P}_h, \ \int_e[v]_e\diff s = 0 \ \forall e\in\mc{E}_h\right\},
\end{eqnarray*}
respectively. Now we define the local and global Kouhia-Stenberg type virtual element spaces as follows:
\begin{dis}
\vV(K) = V_{nc}(K) \times V_{c}(K), \qquad \vV_h = V_{h,nc}\times V_{h,c}.
\end{dis}
The degrees of freedom for $\vV(K)$ can be chosen as, for $\vv\in\vV(K)$ with $\vv = (v_1,v_2)$,
\begin{eqnarray}
& \bullet & \tr{the moments $\frac{1}{|e|}\int_ev_1\diff s$ for each edge $e$ of $K$}, \label{eqn:KSLDOF1} \\
& \bullet & \tr{the values of $v_2$ at each vertex of $K$}, \label{eqn:KSLDOF2}
\end{eqnarray}
The degrees of freedom for $\vV_h$ can be chosen as, for $\vv\in\vV_h$ with $\vv = (v_1,v_2)$,
\begin{eqnarray}
& \bullet & \tr{the moments $\frac{1}{|e|}\int_ev_1\diff s$ for each interior edge $e$}, \label{eqn:KSGDOF1}\\
& \bullet & \tr{the values of $v_2$ at each interior vertex}. \label{eqn:KSGDOF2}
\end{eqnarray}
Given $\vv\in (H_0^1(\Omega)\cap H^2(\Omega))^2$, we denote by $I_h\vv$ the global interpolant of $\vv$, that is, $I_h\vv$ is defined by the unique function in $\vV_h$ such that $\chi_i(\vv - I_h\vv) = 0$ for any $i = 1,2,\cdots,\dim\vV_h$, where $\chi_i$ is the operator that associates the $i$-th degree of freedom of $\vV_h$. Then the following result holds \cite{beirao2013basic,MR3507277}.

\begin{lemma}\label{lem:KSInterpError}
Let $I_h$ be the interpolation operator as defined above. There exists a positive constant $C$ independent of $h$ such that for any $\vv\in (H_0^1(\Omega)\cap H^2(\Omega))^2$ and any $K\in\mc{P}_h$,
\begin{dis}
\|\vv - I_h\vv\|_{0,K} + h_K|\vv - I_h\vv|_{1,K} \leq Ch_K^2|\vv|_{2,K}.
\end{dis}
\end{lemma}

\subsection{Discrete problem}\label{subsec:KSDiscreteProb}

We first define a local projection operator on each element in $\mc{P}_h$. Let $K\in\mc{P}_h$. We define $\Pi_h^K:\vV(K)\to (\mb{P}_1(K))^2$ as the solution of 
\begin{dis}
\left\{\begin{array}{l}
\int_K\vvep(\Pi_{h}^K\vv):\vvep(\vq)\diff\vx = \int_K\vvep(\vv):\vvep(\vq)\diff\vx \quad \forall \vq\in (\mb{P}_1(K))^2 \\
\int_K\rot\Pi_{h}^K\vv\diff\vx = \int_K\rot\vv\diff\vx \\
\int_{\pd K}\Pi_{h}^K\vv\diff s = \int_{\pd K}\vv\diff s
\end{array}\right.
\end{dis}
Note that $\Pi_h^K\vv$ is computable for any $\vv\in\vV(K)$ from the local degrees of freedom \eqref{eqn:KSLDOF1}-\eqref{eqn:KSLDOF2}. Let $\Pi_0^{K}:L^2(K)\to \mb{P}_0(K)$ be the $L^2$-projection operator. Let $S^K:\vV(K)\times\vV(K)\to\R$ be a bilinear form such that
\begin{dis}
S^K(\vu,\vv) = \sum_{i=1}^{\dim\vV(K)}\chi_i(\vu)\chi_i(\vv),
\end{dis}
where $\chi_i$ is the operator associated with the $i$-th local degrees of freedom of $\vV(K)$. We then define the local bilinear form $a_h^K:\vV(K)\times\vV(K)\to\R$ by
\begin{eqnarray*}
a_h^K(\vu,\vv) & = & a_{\mu,h}^K(\vu,\vv) + a_{\lambda,h}^K(\vu,\vv), \quad \tr{where} \\
a_{\mu,h}^K(\vu,\vv) & = & 2\mu\int_K\vvep(\Pi_h^K\vu):\vvep(\Pi_h^K\vv)\diff\vx+ S^K(\vu - \Pi_h^K\vu, \vv - \Pi_h^K\vv), \\
a_{\lambda,h}^K(\vu,\vv) & = & \lambda\int_K(\Pi_0^{K}\div\vu)(\Pi_0^{K}\div\vv)\diff\vx.
\end{eqnarray*}
Following the arguments in \cite{beirao2013basic}, it is easy to see that the bilinear form $a_h^K(\cdot,\cdot)$ satisfies the consistency and the stability:
\begin{itemize}
\item (Consistency) $a_h^K(\vp,\vv) = a^K(\vp,\vv)$ for any $\vv\in\vV(K)$ and $\vp\in(\mb{P}_1(K))^2$;
\item (Stability) there exist two positive constant $c_*$ and $c^*$, independent of $h$ and of $K$, such that $c_*a_{\mu}^K(\vv,\vv) \leq a_{\mu,h}^K(\vv,\vv) \leq c^*a_{\mu}^K(\vv,\vv)$ for any $\vv\in\vV(K)$.
\end{itemize}
Now the global discrete bilinear form $a_h:\vV_h\times\vV_h\to\R$ is defined by
\begin{dis}
a_h(\vu_h,\vv_h) = \sum_{K\in\mc{P}_h}a_h^K(\vu_h,\vv_h), \quad \forall \vu_h,\vv_h\in\vV_h.
\end{dis}

We next construct the discrete loading term. We define $\vf_h$ on each element $K\in\mc{P}_h$ as the $(L^2(K))^2$-projection of $\vf$ on the space of piecewise constant, that is,
\begin{dis}
\vf_h|_K = \frac{1}{|K|}\int_K\vf\diff\vx, \quad \forall K\in\mc{P}_h,
\end{dis}
where $|K|$ denotes the area of $K$. We then define the discrete loading term $\inn{\vf_h,\cdot}$ as follows:
\begin{dis}
\inn{\vf_h,\vv_h} = \sum_{K\in\mc{P}_h}\int_K\vf_h\cdot\what{\vv}_h\diff\vx, \quad \forall \vv_h\in\vV_h,
\end{dis}
where $\hat{\vv}_h$ is defined by
\begin{dis}
\what{\vv}_h|_K = \left(\frac{1}{N_{K}}\sum_{\substack{e\in\mc{E}_h\\ e\sus \pd K}}\frac{1}{|e|}\int_ev_{h,1}\diff s,\ \frac{1}{N_{K}}\sum_{\substack{\vx\in \mc{V}_h\\ \vx\in K}}v_{h,2}(\vx)\right), \quad \forall K\in\mc{P}_h
\end{dis}
with $\vv_h = (v_{h,1},v_{h,2})$. Then the following lemma for the approximation of the loading term $(\vf,\cdot)_{0,\Omega}$ can be found in \cite{beirao2013basic,MR3507277}.

\begin{lemma}\label{lem:KSloading}
Suppose that $\vf\in (L^2(\Omega))^2$. Then there exists a positive constant $C$ independent of $h$ such that
\begin{dis}
\left|\inn{\vf_h,\vv_h} - (\vf,\vv_h)_{0,\Omega}\right| \leq Ch\|\vf\|_{0,\Omega}|\vv_h|_{1,h} \quad \forall \vv_h \in\vV_h.
\end{dis}
\end{lemma}

With the above preparations, we state the following virtual element discretization of the problem \eqref{eqn:ModelProbVar}: Find $\vu_h\in\vV_h$ such that
\begin{eqn}\label{eqn:DiscreteProbKS}
a_h(\vu_h,\vv_h) = \inn{\vf_h,\vv_h} \quad \forall \vv_h\in\vV_h.
\end{eqn}

\subsection{Preliminary results}

In order to study the error estimate for the proposed method \eqref{eqn:DiscreteProbKS}, we first need to prove a discrete inf-sup condition: there exists a positive constant $\beta$ independent of $h$ such that
\begin{dis}
\inf_{q_h\in Q_h}\sup_{\vv_h\in\vV_h}\frac{(\div_h\vv_h,q_h)_{0,\Omega}}{\|q_h\|_{0,\Omega}|\vv_h|_{1,h}} \geq \beta,
\end{dis}
where $Q_h$ is defined by
\begin{dis}
Q_h = \left\{q_h\in L_0^2(\Omega) : q_h|_{K}\in\mb{P}_0(K) \ \forall K\in\mc{P}_h\right\}.
\end{dis}

To show this, we use the macroelement technique \cite{MR1076437,MR725982,MR2609344}. Here we follow the result in \cite[Section 4]{MR2609344}, with slight modification.

We first summarize some definitions and notations. A \emph{macroelement} $M$ is a connected collection of polygonal elements satisfying the regularity assumptions (ii) and (iii) in \Cref{assum:MeshReg}. For a macroelement $M$, we denote by $\mc{E}_M^i$ and $\mc{V}_M^i$ the set of all interior edges and vertices in $M$, respectively. 

A collection of macroelements $\mc{M}_h$ is called a \emph{macroelement partition} of $\mc{P}_h$ if every element is contained in at least one macroelement in $\mc{M}_h$, that is, for each $K\in\mc{P}_h$ there exists $M\in\mc{M}_h$ such that $K\sus M$. 

A macroelement $M$ is said to be \emph{equivalent} to a reference macroelement $\widehat{M}$ if there exists a continuous bijection $F_M:\widehat{M}\to M$ such that the following are true:
\begin{itemize}
\item $F_M(\widehat{M}) = M$.
\item If $\widehat{M}$ consists of elements $\widehat{K}_1,\cdots,\widehat{K}_m$, then $M$ consists of elements $K_1,\cdots,K_m$ such that $F_M(\widehat{K}_i) = K_i$ and both $\widehat{K}_i$ and $K_i$ have the same number of edges.
\item $F_M|_{\widehat{K}_j} = F_{K_j}\circ F_{\widehat{K}_j}^{-1}$ for each $j = 1,2,\cdots,m$.
\item Both $\widehat{M}$ and $M$ have the same number of interior/boundary edges.
\item Both $\widehat{M}$ and $M$ have the same number of interior/boundary vertices.
\end{itemize}
We say that two macroelements are equivalent if they are equivalent to the same reference macroelement.  

Given macroelement $M$, we define local spaces $\vV(M)$ and $Q(M)$ as follows:
\begin{eqnarray*}
V_{nc}(M) & = & \left\{v\in L^2(M) : v|_K\in V_{nc}(K) \ \forall K\sus M, \ \int_e[v]_e\diff s = 0 \ \forall e\in\mc{E}_M\right\}, \\
V_{c}(M) & = & \left\{v\in H_0^1(M) : v|_K\in V_{c}(K) \ \forall K\sus M\right\}, \\
\vV(M) & = & V_{nc}(M)\times V_{c}(M), \\
Q(M) & = & \left\{q\in L^2(M) : q|_K\in\mb{P}_0(K) \ \forall K\sus M\right\},
\end{eqnarray*}
where $\mc{E}_M$ denotes the set of all edges in the macroelement $M$. 

Under the definitions above, we state the macroelement condition as follows. Note that it is essentially identical to Theorem 4.1 of \cite{MR2609344}, and therefore we skip the proof of the theorem here.

\begin{theorem}[Macroelement condition]\label{thm:MacroCond}
Let $\mc{M}_h$ be a macroelement partition of $\mc{P}_h$. Suppose that there exists a fixed set of equivalent classes $\Sigma_1,\cdots,\Sigma_l$ of the macroelements and a positive integer $L$ indpendent of $h$ such that
\begin{enumerate}[label=(\roman*)]
\item For each $M\in\Sigma_i$, $i = 1,\cdots,l$, the space
\begin{eqn}\label{eqn:DivNullDef}
N(M) := \left\{q\in Q(M) : \sum_{K\sus M}\int_K\div\vv q\diff\vx = 0 \ \forall \vv\in \vV(M)\right\}.
\end{eqn}
is one-dimensional, consisting of functions that are constant on $M$.
\item For each $M\in\mc{M}_h$ there exists $i\in\{1,\cdots,l\}$ such that $M\in\Sigma_i$.
\item Each $K\in\mc{P}_h$ is contained in at least one and not more than $L$ macroelements of $\mc{M}_h$. 
\end{enumerate}
Then there exists a positive constant $\beta$ independent of $h$ such that
\begin{dis}
\sup_{\substack{\vv_h\in\vV_h \\ \vv_h\neq\zz}}\frac{(\div_h\vv_h,q_h)_{0,\Omega}}{|\vv_h|_{1,h}} \geq \beta\|q_h\|_{0,\Omega} \quad \forall q_h\in Q_h.
\end{dis}
\end{theorem}

We next define a macroelement partitioning $\mc{M}_h$ and equivalent classes $\Sigma_1,\cdots\Sigma_l$ of macroelements satisfying the assumptions in \Cref{thm:MacroCond}. For each interior vertex $\vx$ in $\mc{P}_h$, let $M_{\vx}$ be the macroelement consisting of elements $K_1(\vx),\cdots,K_{k(\vx)}(\vx)$ having $\vx$ as a common vertex, ordered counterclockwise about the vertex $\vx$ (see, for example, \Cref{fig:Macroelement}). We then define
\begin{dis}
\mc{M}_h := \{M_{\vx} : \vx\in\mc{V}_h^i\}.
\end{dis}
Then $\mc{M}_h$ is clearly a macroelement partitioning of $\mc{P}_h$ satisfying the condition (iii) in \Cref{thm:MacroCond}. 


\begin{figure}
\begin{center}
\begin{tikzpicture}[scale = 0.75]
\fill (0,0) circle [radius = 0.075];
\draw [thick] (0,0) -- (1.5,0.75);
\draw [thick] (0,0) -- (-1,-1.5) -- (-2,-2);
\draw [thick] (0,0) -- (-2,-0.5) -- (-2,-2);
\draw [thick] (0,0) -- (-1,0.75);
\draw [thick] (0,0) -- (1,-1) -- (1,-1.5);
\draw [thick] (-1,0.75) -- (-1.5,1.5) -- (-3,0.5) -- (-2,-0.5);
\draw [thick] (-1.5,1.5) -- (-1,2) -- (0.5,2) -- (1.5,1.5) -- (1.5,0.75);
\draw [thick] (1.5,0.75) -- (2.5,0) -- (2,-0.5) -- (1,-1);
\draw [thick] (1,-1.5) -- (0,-2) -- (-1,-1.5);
\node [above] at (0,0) {$\vx$};
\node at (0,1) {$K_1$};
\node at (-1,0.25) {$K_2$};
\node at (-1.25,-0.75) {$K_3$};
\node at (0,-1) {$K_4$};
\node at (0.75,-0.125) {$K_5$};
\end{tikzpicture}
\caption{}
\label{fig:Macroelement}
\end{center}
\end{figure}

Consider two macroelements $M_{\vx}$ and $M_{\vx'}$ in $\mc{M}_h$. Two macroelements are clearly equivalent if the following are true:
\begin{enumerate}[label=(\roman*)]
\item $k(\vx) = k(\vx')$, that is, the number of polygons in $M_{\vx}$ is equal to the number of polygons in $M_{\vx'}$. 
\item $M_{\vx}$ and $M_{\vx'}$ have the same number of interior/boundary edges and vertices.
\item For $i = 1,2,\cdots,k(\vx)$ $( = k(\vx'))$, the polygons $K_i(\vx)$ and $K_i(\vx')$ have the same number of edges and vertices.
\item The number of edges in $\pd K_i(\vx) \cap \pd K_{i+1}(\vx)$ is equal to the number of edges in $\pd K_i(\vx')\cap \pd K_{i+1}(\vx')$, for $i = 1,2,\cdots,k(\vx)$ modulo $k(\vx)$. 
\end{enumerate}
As mentioned in \Cref{subsec:MeshReg}, since the minimum angle of the triangular decomposition of any polygon in $\mc{P}_h$ is uniformly bounded below by a positive constant controlled by $\rho$, there exists $\ol{k}\in\N$ independent of $h$ such that $k(\vx) \leq \ol{k}$ for any $\vx\in\mc{V}_h^i$. Moreover, since each polygon in $\mc{P}_h$ has at most $N$ edges and vertices, where $N\in\N$ is independent of $h$, each macroelement $M_{\vx}$ has at most $\ol{k}N$ edges and vertices. Therefore there exists at most $\ell$ equivalent class $\Sigma_1,\cdots,\Sigma_{\ell}$ in $\mc{M}_h$, where $\ell$ is a positive integer depending only on $\ol{k}$ and $N$. Thus the condition (ii) in \Cref{thm:MacroCond} is true. Now it remains to show that the condition (i) in \Cref{thm:MacroCond} is also true.

\begin{lemma}\label{lem:MacroPartitionCond1}
Consider the macroelement $M = M_{\vx}$, for $\vx\in\mc{V}_h^i$. Then the space $N_M$ in \eqref{eqn:DivNullDef} is one-dimensional, consisting of functions that are constant on $M$.
\end{lemma}

\begin{proof}
We follow the arguement in the proof of \cite[Lemma 4.3]{MR1343077}. Let $\ol{\vx}\in\mc{V}_h^i$ be fixed and consider the macroelement $M_{\ol{\vx}}\in\mc{M}_h$ consisting of polygons $K_1,\cdots,K_k$, with $k = k(\ol{\vx})$, ordered counterclockwise about the vertex $\ol{\vx}$. Let $e_1,\cdots,e_k$ be the interior edges in $M$ having $\ol{\vx}$ as a common vertex and satisfying $e_i\sus \pd K_i \cap \pd K_{i+1}$ for each $i = 1,\cdots,k$ modulo $k$. Let $q\in N(M)$. If $\vv\in\vV(M)$ with $\vv = (v_1,v_2)$ satisfies $v_2 \equiv 0$, $\int_{e_i}v_1\diff s = |e_i|$ for each $i$, and $\int_ev_1\diff s = 0$ for any other interior edge $e$ in $M$, then 
\begin{dis}
0 = \sum_{i=1}^k\int_{K_i}\div\vv q\diff\vx = \sum_{i=1}^k\int_{\pd K_i}\vv\cdot\vn q\diff s = n_{e_i,1}|e_i|(q_i - q_{i+1}),
\end{dis}
where $\vn_{e_i} = (n_{e_i,1},n_{e_i,2})$ and $q_i = \frac{1}{|K_i|}\int_{K_i}q\diff\vx$ for each $i = 1,\cdots,k$. Thus $q_i = q_{i+1}$ unless $n_{e_i,1} = 0$, for any $i = 1,\cdots,k$ modulo $k$. Note that there exist at most two edges in $\{e_1,\cdots,e_k\}$ such that the normal vectors at the edges are parallel to $y$-axis. 

If there is one index $i\in\{1,\cdots,k\}$ such that $n_{e_i,1} = 0$, then we obtain $q_1 = \cdots = q_i$ and $q_{i+1} = \cdots = q_k = q_1$. Thus $q$ must be constant on $M$. 

If there are two indices $i$ and $i'$ in $\{1,\cdots,k\}$ such that $n_{e_i,1} = 0$ and $n_{e_{i'},1} = 0$ (we may assume that $i < i' < k$), then $q_{i+1} = \cdots = q_{i'}$ and $q_{i'+1} = \cdots = q_k = q_1 = \cdots = q_i$ but $q_i \neq q_{i'}$ in general. In this case, consider $\vv\in\vV(M)$ with $\vv = (v_1,v_2)$ satisfying $v_1 \equiv 0$, $v_2(\ol{\vx}) = 1$, and $v_2(\vx) = 0$ for any other interior vertex $\vx$ in $M$. Then, since $q_{i+1} = \cdots = q_{i'}$ and $q_{i'+1} = \cdots = q_k = q_1 = \cdots = q_i$,
\begin{eqnarray*}
0 & = & \sum_{j=1}^k\int_{K_j}\div\vv q\diff\vx = \sum_{j=1}^k\int_{\pd K_j}\vv\cdot\vn q\diff s \\
& = & \frac{1}{2}n_{e_i,2}|e_i|(q_i - q_{i+1}) + \frac{1}{2}n_{e_{i'},2}|e_{i'}|(q_{i'} - q_{i'+1}) \\
& = & \frac{1}{2}\left(n_{e_i,2}|e_i| - n_{e_{i'},2}|e_{i'}|\right)(q_i - q_{i+1}).
\end{eqnarray*}
Since $\vn_{e_i}$ and $\vn_{e_{i'}}$ are unit vectors such that $n_{e_i,1} = n_{e_{i'},1} = 0$ and $e_1,\cdots,e_k$ are edges having $\ol{\vx}$ as a common vertex and ordered counterclockwise about $\ol{\vx}$, we have $n_{e_i,2} = -n_{e_{i'},2}$ and $n_{e_i,2} \neq 0$. Thus $q_i = q_{i+1}$ and $q$ must be constant on $M$.
\end{proof}

As a corollary, we now obtain that the Kouhia-Stenberg type virtual element satisfies the discrete inf-sup condition. 

\begin{corollary}[Discrete inf-sup condition]\label{cor:InfSup}
There exists a positive constant $\beta$ independent of $h$ such that 
\begin{dis}
\sup_{\substack{\vv_h\in\vV_h \\ \vv_h\neq\zz}}\frac{(\div_h\vv_h,q_h)_{0,\Omega}}{|\vv_h|_{1,h}} \geq \beta\|q_h\|_{0,\Omega} \quad \forall q_h\in Q_h.
\end{dis}
\end{corollary}

Using \Cref{cor:InfSup} and the classical arguments (see, for example, Proposition 2.5 in Chapter 2 of \cite{MR1115205}), we have the following property.

\begin{corollary}\label{cor:KSdivInterpol}
Let $\vu\in (H_0^1(\Omega)\cap H^2(\Omega))^2$. Then there exists $\vu_I\in\vV_h$ such that
\begin{eqn}\label{eqn:KSNewInterpolant}
(\div_h\vu_I,q_h)_{0,\Omega} = (\div\vu,q_h)_{0,\Omega} \ \forall q_h\in Q_h, \quad |\vu - \vu_I|_{1,h} \leq Ch|\vu|_{2,\Omega},
\end{eqn}
where $C$ is a positive constant independent of $h$. 
\end{corollary}

\begin{proof}
We follow the argument in the proof of Proposition 2.5 in Chapter 2 of \cite{MR1115205}. Let $\vv_h = I_h\vu$, where $I_h$ is the interpolation operator in \Cref{lem:KSInterpError}. Then, from \Cref{cor:InfSup}, there exists $\vr_h\in\vV_h$ such that
\begin{dis}
\left(\div_h\vr_h,q_h\right)_{0,\Omega} = \left(\div\vu - \div_h\vv_h,q_h\right)_{0,\Omega} \quad \forall q_h\in Q_h
\end{dis}
and $|\vr_h|_{1,h} \leq C|\vu - \vv_h|_{1,h}$, where $C$ is a positive constant independent of $h$. Define $\vu_I := \vr_h + \vv_h$. Then $\vu_I$ satisfies $(\div_h\vu_I,q_h)_{0,\Omega} = (\div\vu,q_h)_{0,\Omega}$ for any $q_h\in Q_h$. From \Cref{lem:KSInterpError},
\begin{dis}
|\vu - \vu_I|_{1,h} \leq |\vu - \vv_h|_{1,h} + |\vr_h|_{1,h} \leq C|\vu - \vv_h|_{1,h} \leq Ch|\vu|_{2,\Omega}.
\end{dis}
This completes the proof.
\end{proof}

We next prove a discrete version of Korn's inequality. 

\begin{theorem}\label{thm:KornKS}
There exists a positive constant $C$ independent of $h$ such that
\begin{dis}
|\vv_h|_{1,h} \leq C\|\vvep_h(\vv_h)\|_{0,\Omega} \quad \forall \vv_h\in\vV_h.
\end{dis}
\end{theorem}

\begin{proof}
We first define some finite element spaces on the triangulation $\mc{T}_h$ as follows:
\begin{eqnarray*}
M_{h} & = & \left\{w\in L^2(\Omega) : w|_T\in \mb{P}_1(T) \ \forall T\in\mc{T}_h\right\}, \\
X_{h,c} & = & M_h\cap H_0^1(\Omega), \\
X_{h,nc} & = & \left\{w\in M_{h} : \parbox[c]{95mm}{\centering$w$ is continuous at the midpoints of interior edges of $\mc{T}_h$ \\ and $w = 0$ at the midpoints of boundary edges of $\mc{T}_h$} \right\}, \\
\vX_h & = & X_{h,nc}\times X_{h,c}.
\end{eqnarray*}
According to \cite[Lemma 4.5]{MR1343077}, there is a positive constant $C$ independent of $h$ such that
\begin{eqn}\label{eqn:KSKorn001}
|\vw_h|_{H^1(\Omega;\mc{T}_h)} \leq C\|\vvep_h(\vw_h)\|_{0,\Omega} \quad \forall \vw_h\in\vX_h.
\end{eqn}

We next define a subspace $\vW_h$ of $\vX_h$ whose degrees of freedom can be chosen as the same with the degrees of freedom of $\vV_h$. Let $K\in\mc{P}_h$. We define some auxiliary spaces as follows:
\begin{align*}
& \wtilde{W}(K) = \left\{v\in H^1(K) : v|_T\in \mb{P}_1(T) \ \forall T\in\mc{T}^K\right\}, \quad \wtilde{W}_0(K) = H_0^1(K) \cap W(K), \\
& B_D(\pd K) = \left\{g\in C^0(\pd K) : g|_e\in\mb{P}_1(e)\ \forall e\sus\pd K\right\}, \\
& B_N(\pd K) = \left\{g\in L^2(\pd K) : g|_e\in \mb{P}_0(e) \ \forall e\sus\pd K\right\}.
\end{align*}
The local spaces $W_c(K)$ and $W_{nc}(K)$ are defined as follows:
\begin{eqnarray*}
W_c(K) & = & \left\{v\in \wtilde{W}(K)\cap B_D(\pd K) : (\nabla v,\nabla w)_{0,K} = 0 \ \forall w\in \wtilde{W}_0(K)\right\} \\
W_{nc}(K) & = & \left\{v\in \wtilde{W}(K) : \exists g\in B_N(\pd K) \ \tr{such that} \ (\nabla v,\nabla w)_{0,K} = {\textstyle\int_{\pd K}}gw\diff s \ \forall w\in \wtilde{W}_0(K)\right\}.
\end{eqnarray*}
That is, $W_c(K)$ consists of P1-conforming finite element approximate solutions of the Dirichlet problem
\begin{dis}
-\Delta v = 0 \quad \tr{in} \ K, \qquad v = g \quad \tr{on} \ \pd K
\end{dis}
with $g\in B_D(\pd K)$, and $W_{nc}(K)$ consists of P1-conforming finite element approximate solutions of the Neumann problem 
\begin{dis}
-\Delta v = 0 \quad \tr{in} \ K, \qquad \pd v/\pd\vn = g \quad \tr{on} \ \pd K
\end{dis}
with $g\in B_N(\pd K)$. We then define $\vW_h := W_{h,nc}\times W_{h,c}$, where
\begin{eqnarray*}
W_{h,c} & = & \left\{v\in H_0^1(\Omega) : v|_K\in W_c(K) \ \forall K\in\mc{P}_h\right\}, \\
W_{h,nc} & = & \left\{v\in L^2(\Omega) : v|_K\in W_{nc}(K) \ \forall K\in\mc{P}_h, \ \int_e[v]_e\diff s = 0 \ \forall e\in\mc{E}_h\right\}.
\end{eqnarray*}
It is clear that $\vW_h$ is a subspace of $\vX_h$ and its degrees of freedom can be chosen as \eqref{eqn:KSGDOF1}-\eqref{eqn:KSGDOF2}. Let $\Phi: \vV_h\to \vW_h$ be the linear bijection such that both $\vv_h$ and $\Phi(\vv_h)$ have exactly the same values of degrees of freedom. That is, given $\vv_h\in \vV_h$ with $\vv_h = (v_{h,1},v_{h,2})$, we define $\vw_h := \Phi(\vv_h)$ in $\vW_h$ with $\vw_h = (w_{h,1},w_{h,2})$ as
\begin{dis}
\int_ew_{h,1}\diff s = \int_ev_{h,1}\diff s \quad \forall e\in\mc{E}_h, \qquad w_{h,2}(\vx) = v_{h,2}(\vx) \quad \forall \vx\in\mc{V}_h.
\end{dis}
Moreover, there exist two positive constant $c_*$ and $c^*$ independent of $h$ such that
\begin{eqn}\label{eqn:KSKorn002}
c_*|\vv_h|_{1,h} \leq |\Phi(\vv_h)|_{1,h} \leq c^*|\vv_h|_{1,h} \quad \forall \vv_h\in\vV_h.
\end{eqn}
Note that $a_{\mu,h}(\vv_h,\vv_h')$ depends only on the degrees of freedom of $\vv_h,\vv_h'\in\vV_h$ (see \Cref{subsec:KSDiscreteProb}). Thus we can construct a discrete bilinear form $\what{a}_{\mu,h}: \vW_h\times\vW_h\to \R$ satisfying the following properties:
\begin{itemize}
\item If $\vv_h,\vv_h'\in\vV_h$, $\vw_h = \Phi(\vv_h)$ and $\vw_h' = \Phi(\vv_h')$, then 
\begin{eqn}\label{eqn:KSKorn003}
a_{\mu,h}(\vv_h,\vv_h') = \what{a}_{\mu,h}(\vw_h,\vw_h').
\end{eqn}
\item there exist two positive constant $c_*$ and $c^*$ independent of $h$ such that
\begin{eqn}\label{eqn:KSKorn004}
c_*\|\vvep_h(\vw_h)\|_{0,\Omega}^2 \leq \what{a}_{\mu,h}(\vw_h,\vw_h) \leq c^*\|\vvep_h(\vw_h)\|_{0,\Omega}^2 \quad \forall \vw_h\in\vW_h.
\end{eqn}
\end{itemize}
Using \eqref{eqn:KSKorn001}-\eqref{eqn:KSKorn004} and the stability of the bilinear form $a_{\mu,h}(\cdot,\cdot)$, we obtain 
\begin{eqnarray*}
|\vv_h|_{1,h}^2 & \leq & C|\Phi(\vv_h)|_{1,h}^2 \leq C\|\vvep_h(\Phi(\vv_h))\|_{0,\Omega}^2 \\
& \leq & C\what{a}_{\mu,h}(\Phi(\vv_h),\Phi(\vv_h)) = Ca_{\mu,h}(\vv_h,\vv_h) \leq C\|\vvep_h(\vv_h)\|_{0,\Omega}^2
\end{eqnarray*}
for any $\vv_h\in\vV_h$. This conclude the proof of the theorem.
\end{proof}

\subsection{Error analysis}

We are now ready to prove the following convergence theorem.

\begin{theorem}\label{thm:KSVEMconv}
Suppose that $\vf\in (L^2(\Omega))^2$. The discrete problem \eqref{eqn:DiscreteProbKS} has a unique solution $\vu_h\in\vV_h$ and, if $\vu\in (H^2(\Omega)\cap H_0^1(\Omega))^2$ is the solution of \eqref{eqn:ModelProbVar}, then
\begin{eqn}\label{eqn:KSErrorTheorem}
|\vu - \vu_h|_{1,h} \leq Ch\|\vf\|_{0,\Omega},
\end{eqn}
where $C$ is a positive constant independent of $h$ and the Lam\'e constant $\lambda$. 
\end{theorem}

\begin{proof}
The proof is almost the same as the one of \Cref{thm:NCVEMconv}. Since the discrete bilinear form $a_h(\cdot,\cdot)$ is elliptic by \Cref{thm:KornKS} and is bounded, the discrete problem \eqref{eqn:DiscreteProbKS} has a unique solution, say $\vu_h$. Let $\vu_{\pi}$ be the approximation in \Cref{lem:ApproxPoly} and $\vdel_h = \vu_h - \vu_I$, where $\vu_I$ is the function in $\vV_h$ satisfying \eqref{eqn:KSNewInterpolant}. Then, according to \Cref{thm:KornKS} and the consistency of $a_h(\cdot,\cdot)$, 
\begin{eqnarray}
C|\vdel_h|_{1,h}^2 & \leq & \inn{\vf_h,\vdel_h} - (\vf,\vdel_h)_{0,\Omega} - \sum_{K\in\mc{P}_h}\left(a_h^K(\vu_I - \vu_{\pi},\vdel_h) + a^K(\vu_{\pi} - \vu,\vdel_h)\right) \nonumber\\
&& - \sum_{K\in\mc{P}_h}a^K(\vu,\vdel_h) + (\vf,\vdel_h)_{0,\Omega}. \label{eqn:KSError00}
\end{eqnarray}
By \Cref{lem:KSloading}, 
\begin{eqn}
\left|\inn{\vf_h,\vdel_h} - (\vf,\vdel_h)_{0,\Omega}\right| \leq Ch\|\vf\|_{0,\Omega}|\vdel_h|_{1,h}. \label{eqn:KSError01}
\end{eqn}
Note that
\begin{eqnarray*}
&& a_h^K(\vu_I - \vu_{\pi},\vdel_h) + a^K(\vu_{\pi} - \vu,\vdel_h) \\
& = & a_{\mu,h}^K(\vu_I - \vu_{\pi},\vdel_h) + a_{\mu}^K(\vu_{\pi} - \vu,\vdel_h) + a_{\lambda,h}^K(\vu_I - \vu_{\pi},\vdel_h) + a_{\lambda}^K(\vu_{\pi} - \vu,\vdel_h), 
\end{eqnarray*}
From \eqref{eqn:KSNewInterpolant},
\begin{dis}
a_{\lambda,h}^K(\vu_I - \vu_{\pi},\vdel_h) + a_{\lambda}^K(\vu_{\pi} - \vu,\vdel_h) = \lambda(\Pi_0^{K}\div\vu - \div\vu,\div\vdel_h)_{0,K}
\end{dis}
for any $K\in\mc{P}_h$. According to \Cref{lem:ApproxPoly} and \Cref{cor:KSdivInterpol} we obtain
\begin{eqnarray}
&& \left|\sum_{K\in\mc{P}_h}a_h^K(\vu_I - \vu_{\pi},\vdel_h) + a^K(\vu_{\pi} - \vu,\vdel_h)\right| \nonumber\\
& \leq & C\sum_{K\in\mc{P}_h}\left(|\vu_I - \vu_{\pi}|_{1,K} + |\vu - \vu_{\pi}|_{1,K} + \lambda|\Pi_0^{K}\div\vu - \div\vu|_{0,K}\right)|\vdel_h|_{1,K} \nonumber\\
& \leq & Ch\left(|\vu|_{2,\Omega} + \lambda|\div\vu|_{1,\Omega}\right)|\vdel_h|_{1,h}. \label{eqn:KSError02}
\end{eqnarray}
Integrating by parts we obtain
\begin{dis}
\sum_{K\in\mc{P}_h}a^K(\vu,\vdel_h) - (\vf,\vdel_h)_{0,\Omega} = \sum_{K\in\mc{P}_h}\int_{\pd K}\left(\vsig(\vu)\vn_K\right)\cdot\vdel_h\diff s = \sum_{e\in\mc{E}_h}\int_e\vsig(\vu)\vn_e\cdot[\vdel_h]_e\diff s.
\end{dis}
For each $e\in\mc{E}_h$, let $\vP_e^0\vsig(\vu)\vn_e$ be the $L^2$-orthogonal projection of $\vsig(\vu)\vn_e$ onto $(\mb{P}_0(e))^{2}$. Then, since 
\begin{eqn}
\int_e[v_{h,1}]_e\diff s = 0, \ [v_{h,2}]_e = 0, \quad \forall e\in\mc{E}_h, \ \forall \vv_h=(v_{h,1},v_{h,2})\in\vV_h, \label{eqn:KSDiscreteJump}
\end{eqn}
we obtain 
\begin{dis}
\int_e\vsig(\vu)\vn_e\cdot[\vdel_h]_e\diff s = \int_e\left(\vsig(\vu)\vn_e - \vP_e^0\vsig(\vu)\vn_e\right)\cdot[\vdel_h]_e\diff s.
\end{dis}
Let $\vs_e = \vsig(\vu)\vn_e - \vP_e^0\vsig(\vu)\vn_e$, and write $\vs_e = (s_{e,1},s_{e,2})$ and $\vdel_h = (\delta_{h,1},\delta_{h,2})$. Using \eqref{eqn:KSDiscreteJump} again, we obtain
\begin{dis}
\int_e\left(\vsig(\vu)\vn_e - \vP_e^0\vsig(\vu)\vn_e\right)\cdot[\vdel_h]_e\diff s = \int_e s_{e,1}[\delta_{h,1}]_e\diff s = \int_e s_{e,1}[\delta_{h,1} - P_e^0\delta_{h,1}]_e\diff s,
\end{dis}
where $P_e^0\delta_{h,1}$ denotes the $L^2$-orthogonal projection of $\delta_{h,1}$ onto $\mb{P}_0(e)$. From the classical arguement in \cite{MR343661}, if $e\in\mc{E}_h^i$ and $e$ is a common edge of two elements $K^+$ and $K^-$ in $\mc{P}_h$, then
\begin{eqnarray*}
\|\vsig(\vu)\vn_e - \vP_e^0\vsig(\vu)\vn_e\|_{0,e} & \leq & Ch_e^{1/2}\|\vsig(\vu)\|_{1,K^+\cup K^-}, \\
\|[\delta_{h,1} - P_e^0\delta_{h,1}]_e\|_{0,e} & \leq & Ch^{1/2}\left(|\delta_{h,1}|_{1,K^+}^2 + |\delta_{h,1}|_{1,K^-}^2\right)^{1/2}.
\end{eqnarray*}
If $e\in\mc{E}_h^b$, then \eqref{eqn:KSDiscreteJump} implies that $[\delta_{h,1} - P_e^0\delta_{h,1}]_e = 0$. Thus
\begin{eqn}
\left|\sum_{K\in\mc{P}_h}a^K(\vu,\vdel_h) - (\vf,\vdel_h)_{0,\Omega}\right| \leq Ch\left(\|\vu\|_{2,\Omega} + \lambda\|\div\vu\|_{1,\Omega}\right)|\vdel_h|_{1,h}\label{eqn:KSError03}
\end{eqn}
Combining the results \eqref{eqn:KSError00}, \eqref{eqn:KSError01}, \eqref{eqn:KSError02}, and \eqref{eqn:KSError03}, we obtain
\begin{dis}
|\vdel_h|_{1,h}^2 \leq Ch\left(\|\vu\|_{2,\Omega} + \lambda\|\div\vu\|_{1,\Omega} + \|\vf\|_{0,\Omega}\right)|\vdel_h|_{1,h},
\end{dis}
which, together with the regularity estimate \eqref{eqn:regularity}, leads to \eqref{eqn:KSErrorTheorem}. 
\end{proof}

\section{Numerical Experiments}

In this section we present some numerical experiments for the lowest-order nonconforming VEM with stabilizing term introduced in \Cref{sec:NCVEMSTAB} and the Kouhia-Stenberg type VEM introduced in \Cref{sec:KSVEM}. 

Let $\Omega = [0,1]^2$ and $\mu = 1$. Consider the problem \eqref{eqn:ModelProb} where the exact solution is given by
\begin{dis}
\vu(x,y) = \begin{pmatrix}
(\cos(2\pi x) - 1)\sin(2\pi y) + \frac{1}{1+\lambda}\sin(2\pi x)\sin(2\pi y) \\
-(\cos(2\pi y)-1)\sin(2\pi x) + \frac{1}{1+\lambda}x(1-x)y(1-y)
\end{pmatrix}.
\end{dis}

\begin{figure}[!ht]
\centering
\includegraphics[width = 0.25\textwidth]{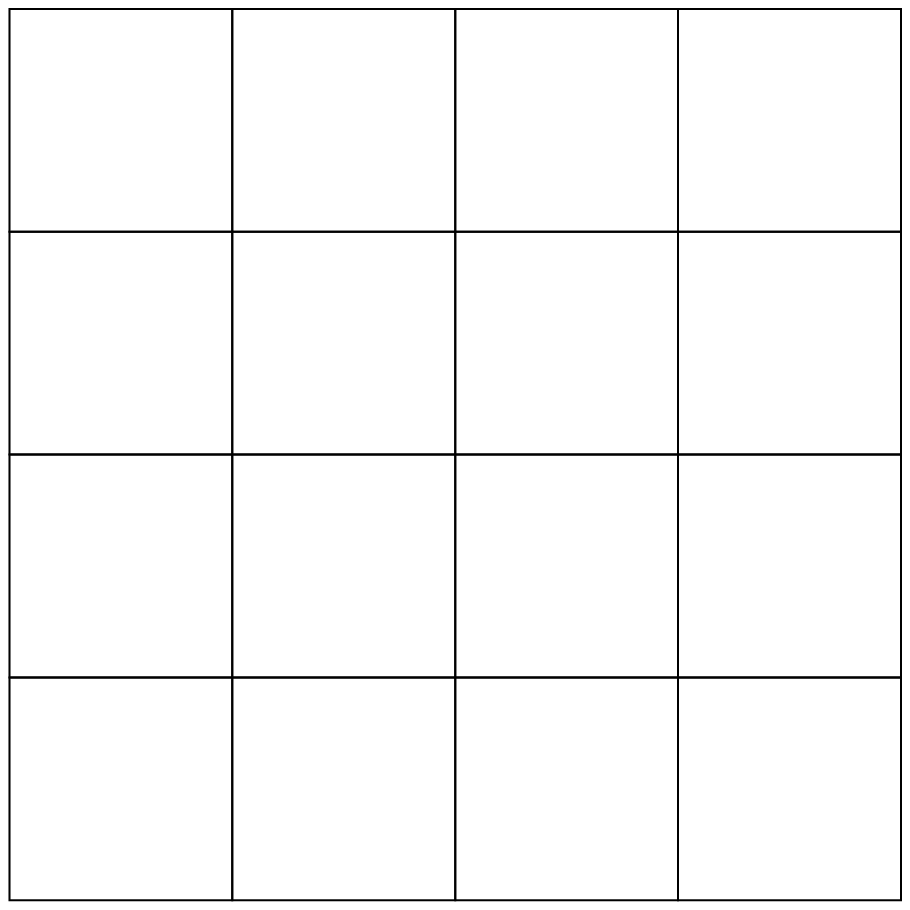}
\includegraphics[width = 0.25\textwidth]{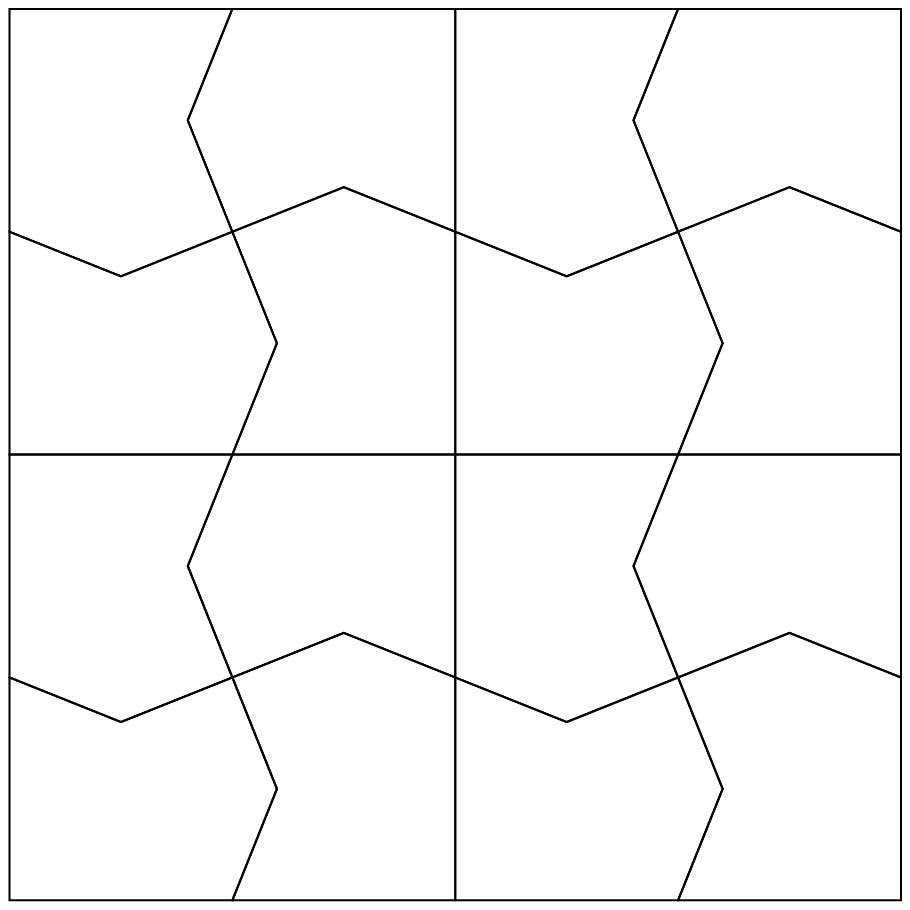}
\includegraphics[width = 0.25\textwidth]{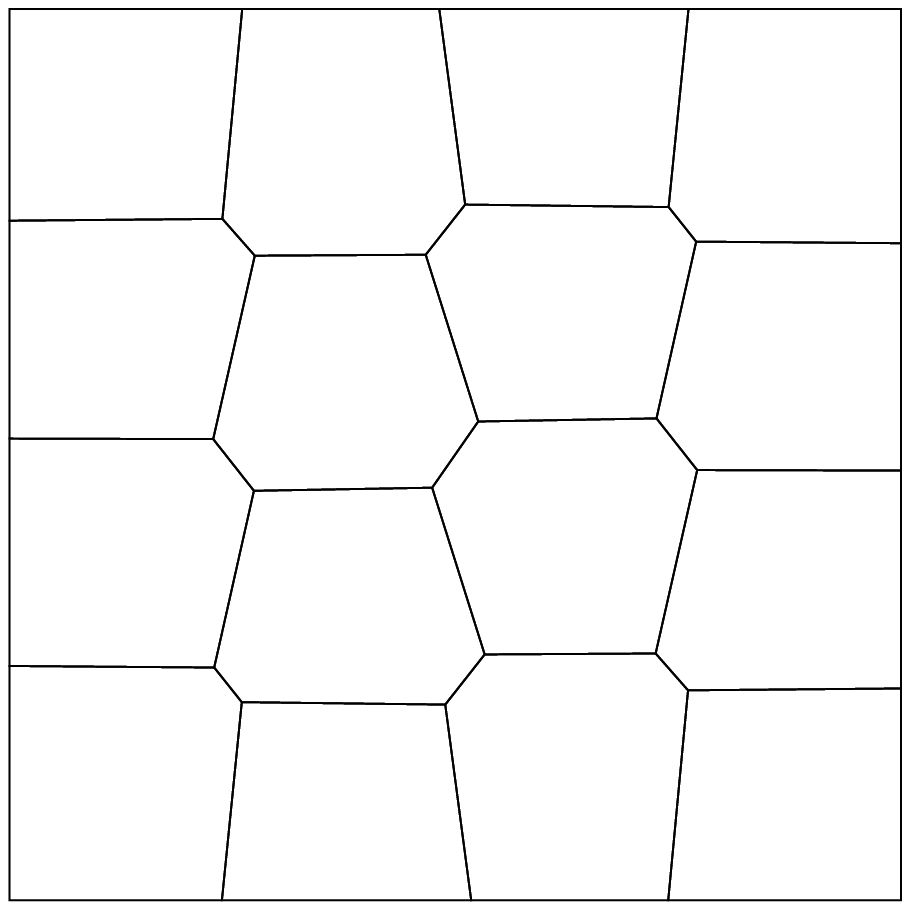}
\includegraphics[width = 0.25\textwidth]{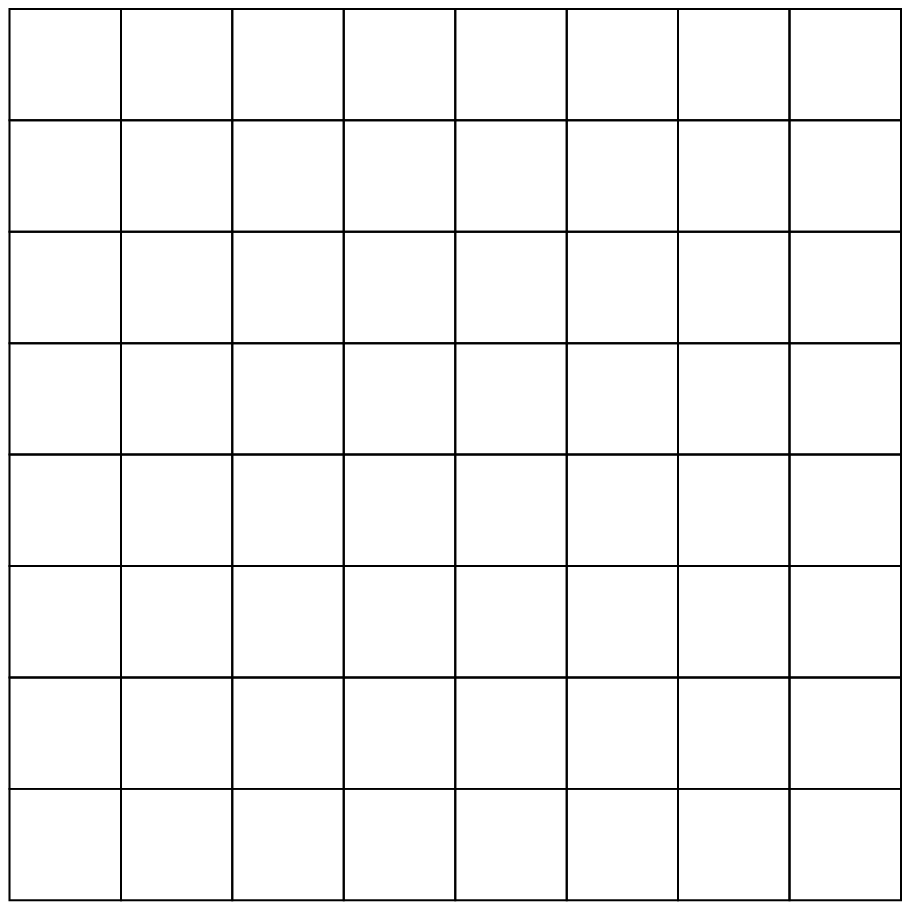}
\includegraphics[width = 0.25\textwidth]{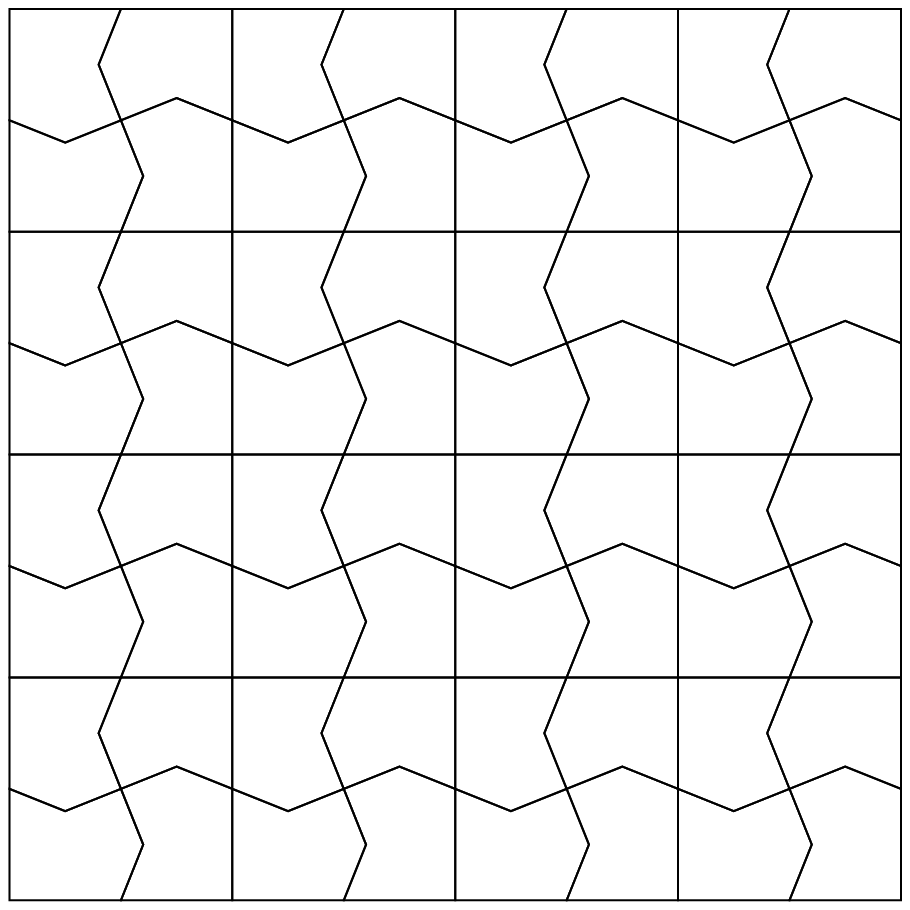}
\includegraphics[width = 0.25\textwidth]{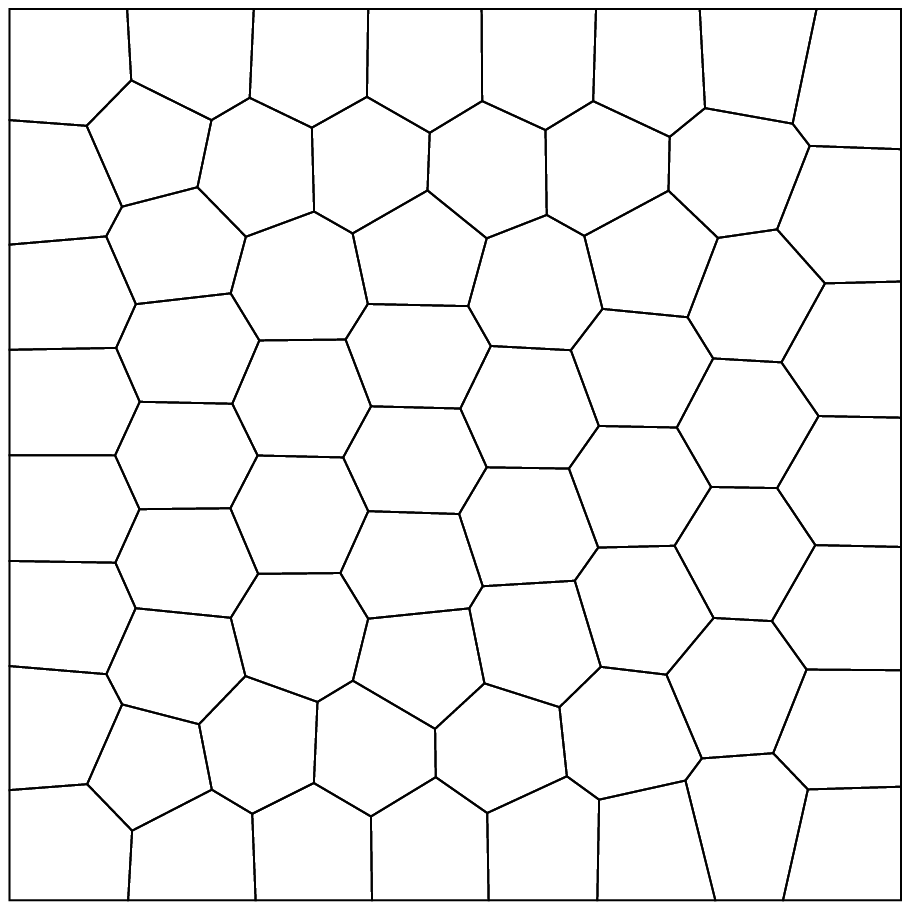}
\caption{The meshes $\mc{P}_h^{(1)}$ (left), $\mc{P}_h^{(2)}$ (middle), and $\mc{P}_h^{(3)}$ (right).}
\label{fig:mesh}
\end{figure}

We consider the following different families of meshes.
\begin{enumerate}[label=(\roman*)]
\item uniform square meshes $\mc{P}_h^{(1)}$ with $h = 1/4, 1/8, 1/16, 1/32, 1/64$,
\item uniform nonconvex hexagonal meshes $\mc{P}_h^{(2)}$ with $h = 1/4, 1/8, 1/16, 1/32, 1/64$,
\item unstructured convex polygonal meshes $\mc{P}_h^{(3)}$ with $h = 1/4, 1/8, 1/16, 1/32, 1/64$.
\end{enumerate}
Some examples of the meshes are shown in \Cref{fig:mesh}. The unstructured convex polygonal meshes are generated from PolyMesher \cite{talischi2012polymesher}. The Lam\'e constant $\lambda$ is taken to be $1$ and $10^4$, respectively.

\subsection{Lowest-order nonconforming element}

\begin{figure}
\centering
\includegraphics[width = 0.4\textwidth]{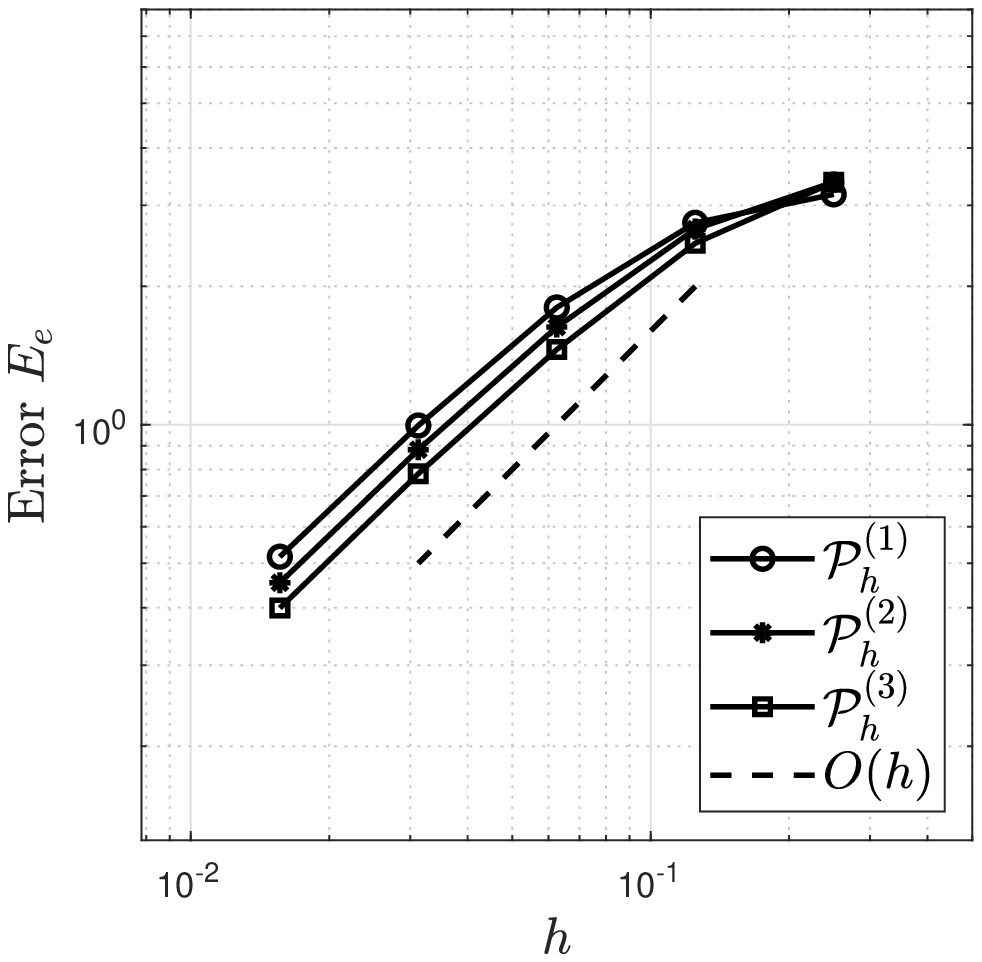}
\includegraphics[width = 0.4\textwidth]{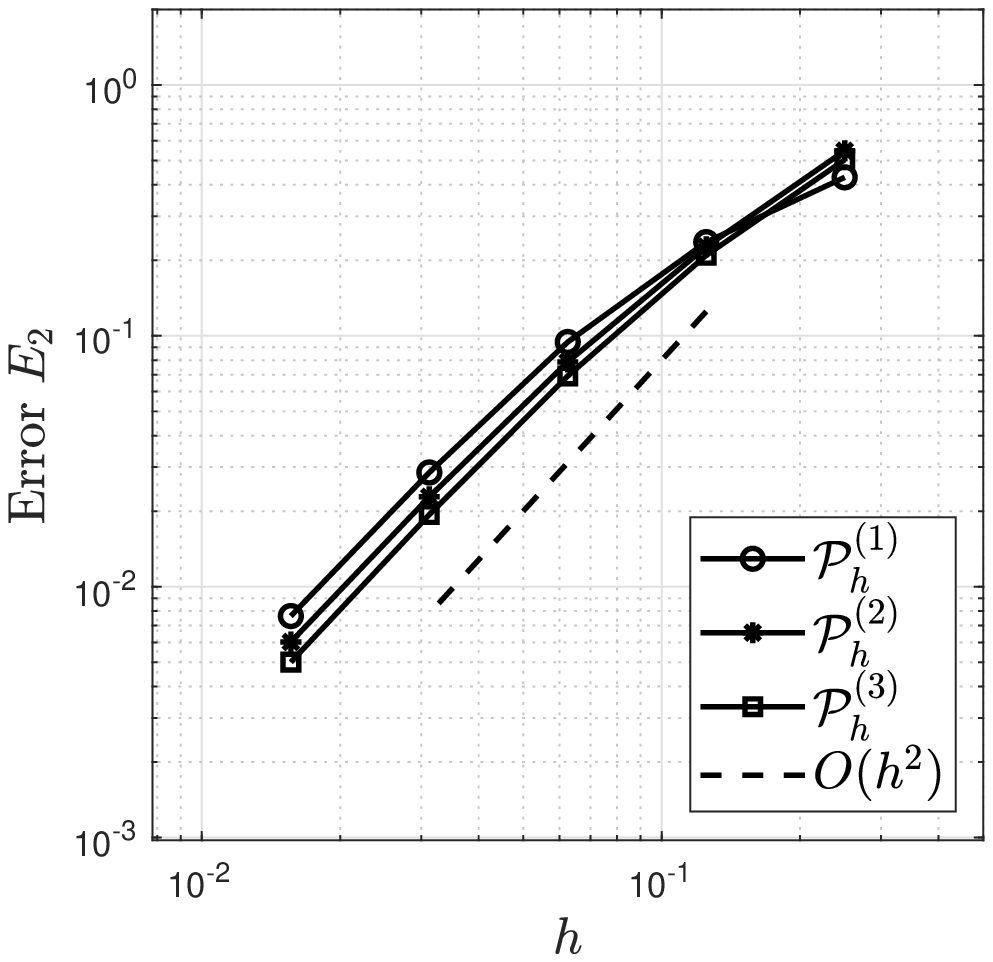}
\caption{The error curves $E_e$ (left) and $E_2$ (right) of test case 1 with $\lambda = 1$.}
\label{fig:errorNC1}
\end{figure}

\begin{figure}
\centering
\includegraphics[width = 0.4\textwidth]{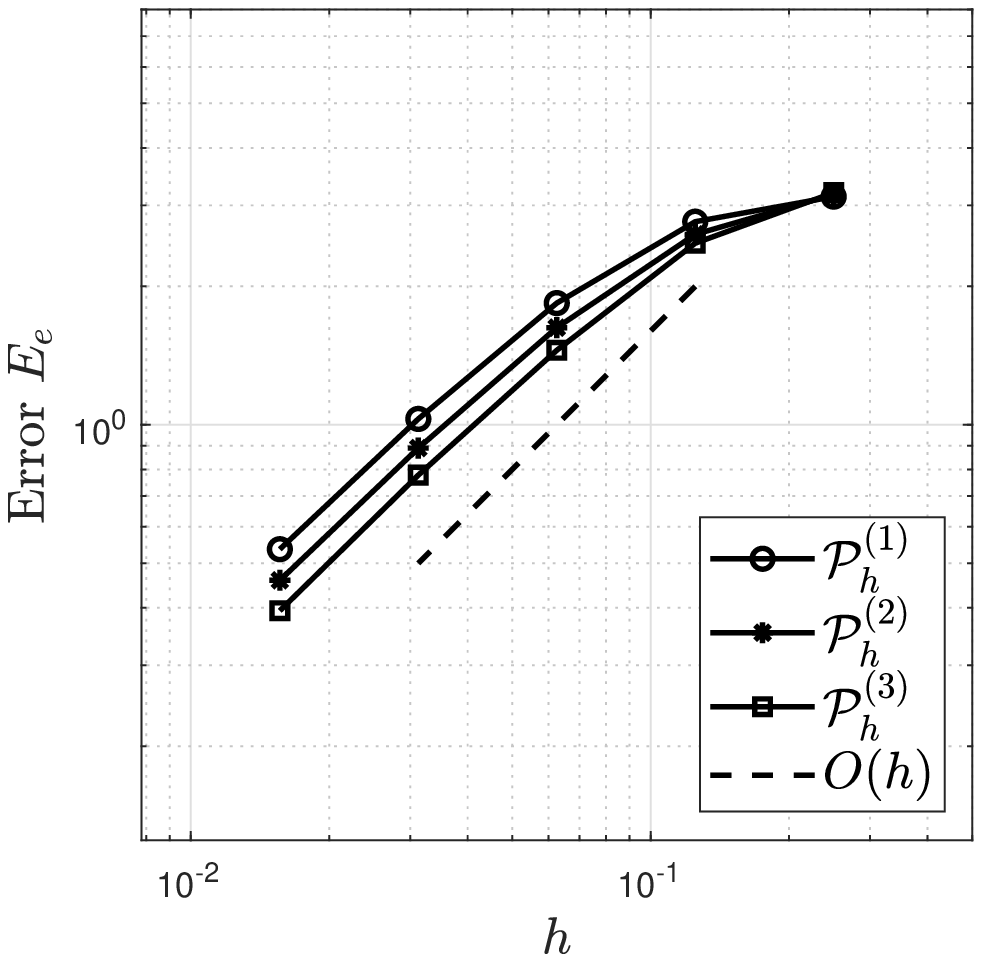}
\includegraphics[width = 0.4\textwidth]{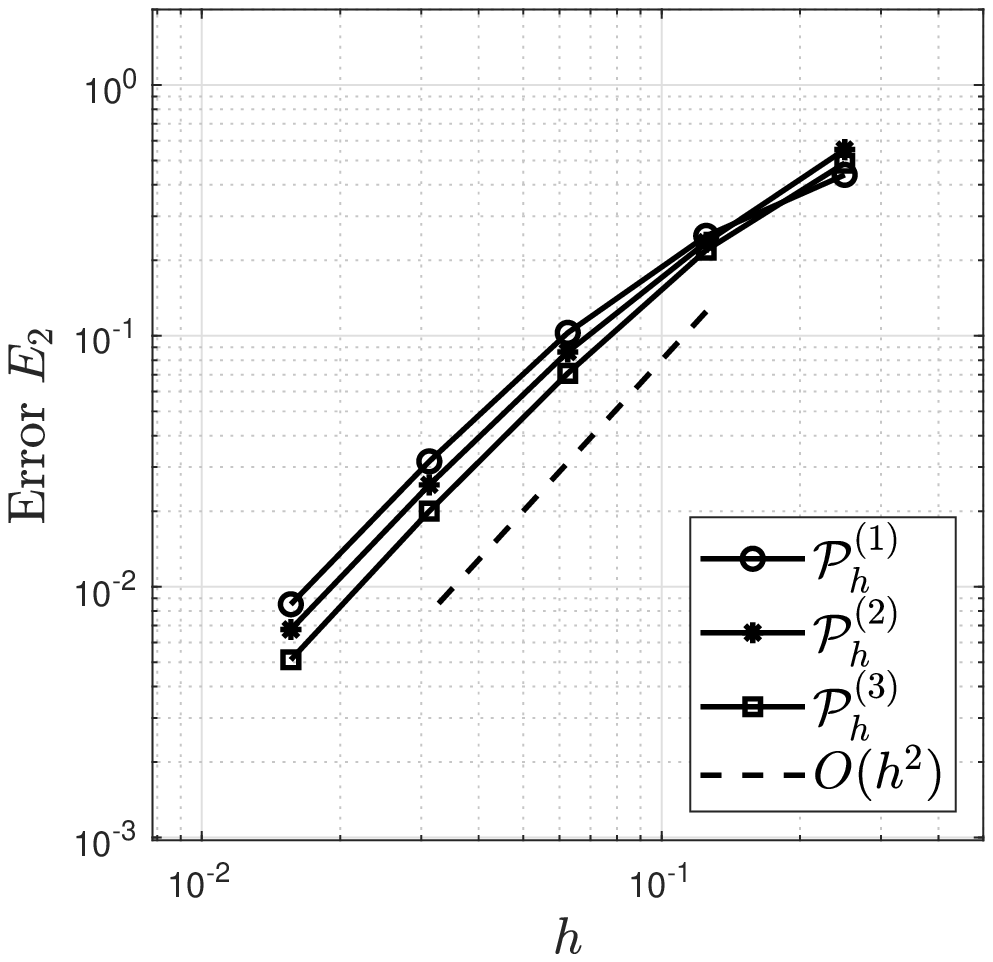}
\caption{The error curves $E_e$ (left) and $E_2$ (right) of test case 1 with $\lambda = 10^4$.}
\label{fig:errorNC2}
\end{figure}

We first implement the method \eqref{eqn:DiscreteProbNCS} with $\gamma = 1$ and compute the errors in the discrete energy norm and $L^2$-norm
\begin{dis}
E_e := a_h(\vu_h - I_h\vu, \vu_h - I_h\vu)^{1/2}, \quad E_2 = \left(\sum_{i=1}^{N_h}h^2|\chi_i(\vu_h) - \chi_i(\vu)|^2\right)^{1/2},
\end{dis}
where $N_h = \dim\vV_h$ and $\chi_i$ is the operator associated with the $i$-th degree of freedom. In \Crefrange{fig:errorNC1}{fig:errorNC2}, we present the error curves versus $h$ for different values of $\lambda$. As shown in these figures, we see that the convergence order of the errors $E_e$ and $E_2$ are $O(h)$ and $O(h^2)$, respectively. Moreover, the convergence order is maintained in the nearly incompressible case ($\lambda = 10^4$). These results are consistent with the convergence rate predicted by the analysis in \Cref{thm:NCVEMconv}. 

\subsection{Kouhia-Stenberg type element}

\begin{figure}
\centering
\includegraphics[width = 0.4\textwidth]{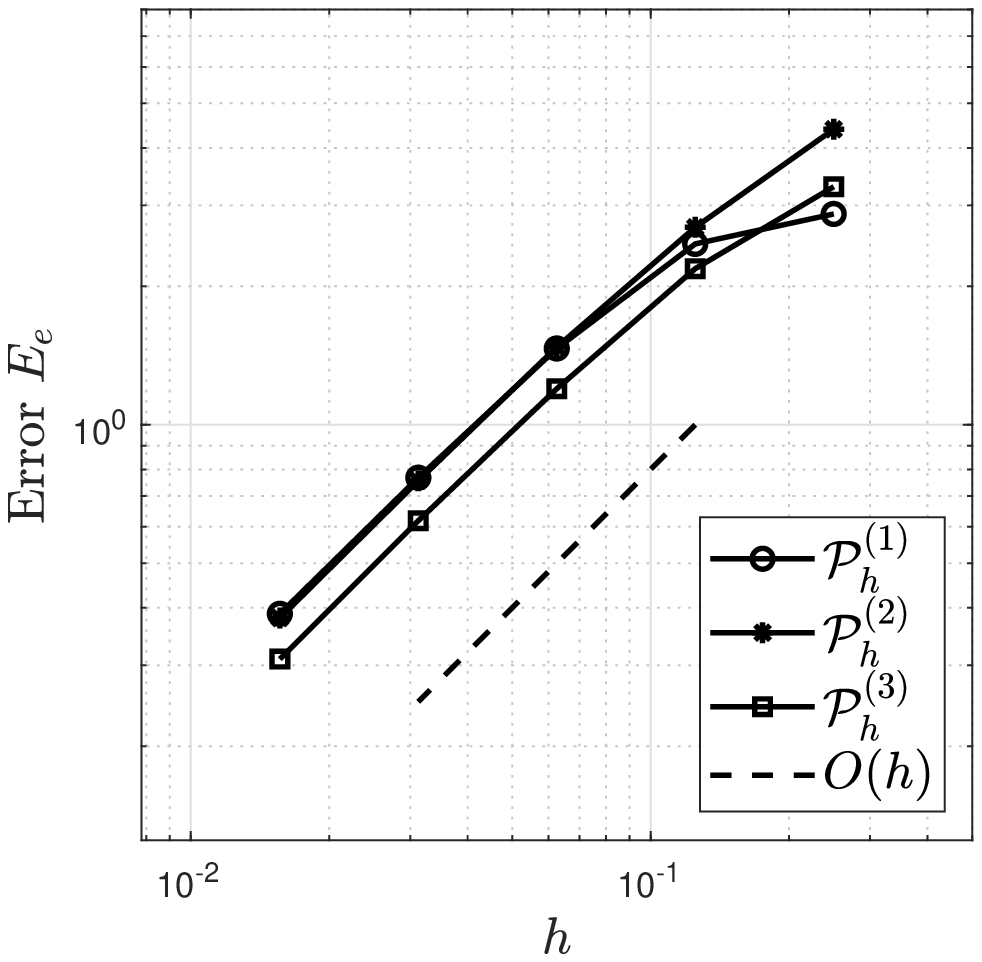}
\includegraphics[width = 0.4\textwidth]{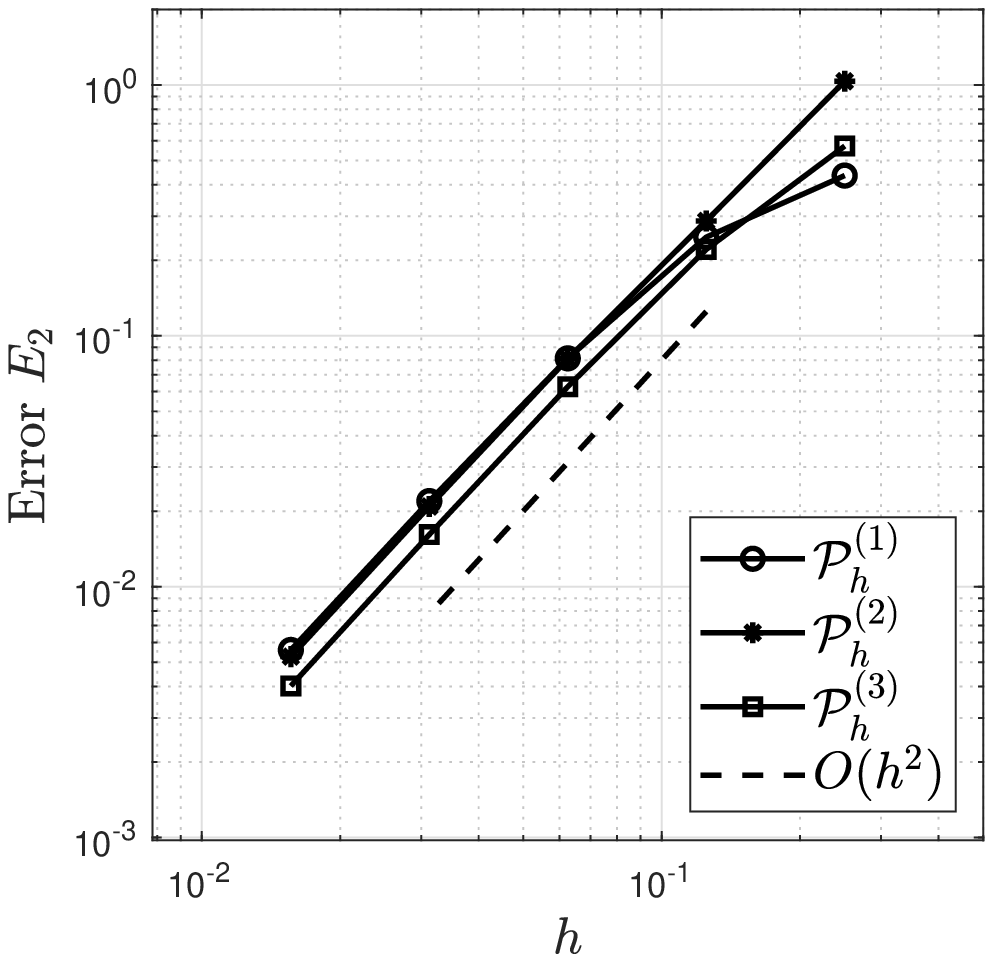}
\caption{The error curves $E_e$ (left) and $E_2$ (right) of test case 2 with $\lambda = 1$.}
\label{fig:errorKS1}
\end{figure}

\begin{figure}
\centering
\includegraphics[width = 0.4\textwidth]{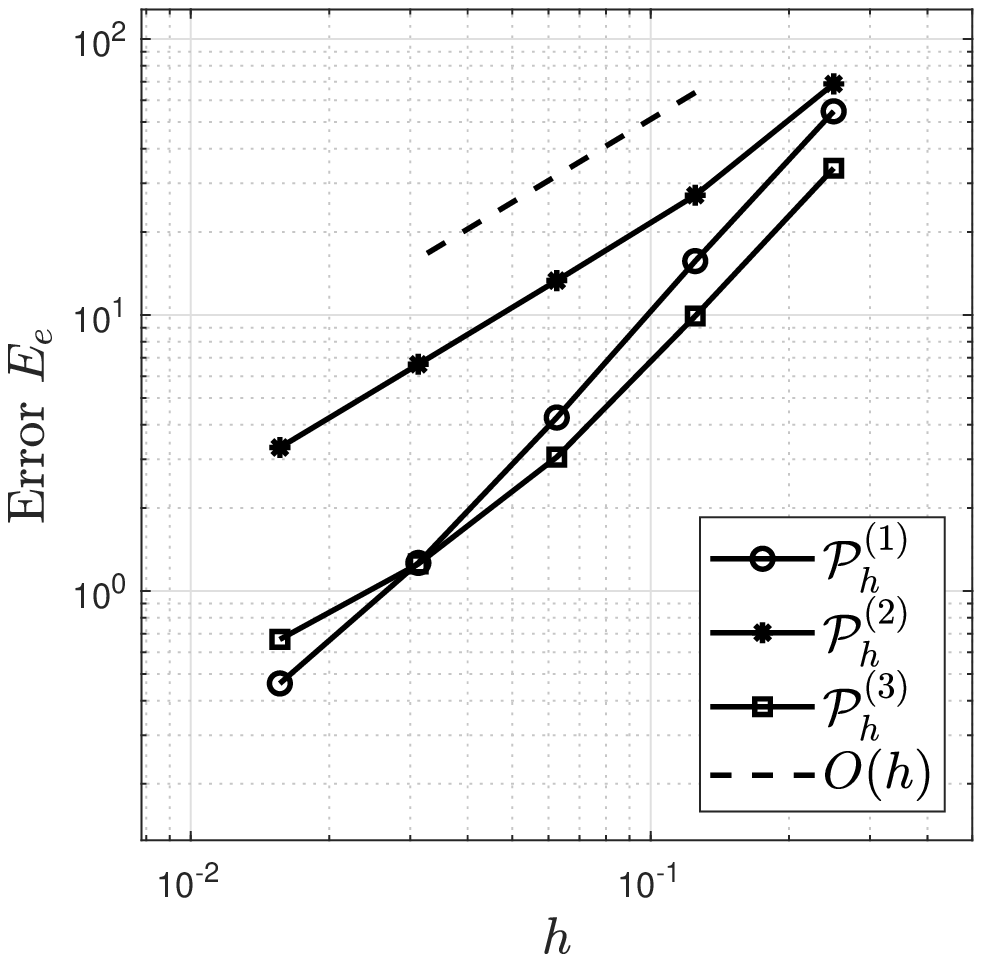}
\includegraphics[width = 0.4\textwidth]{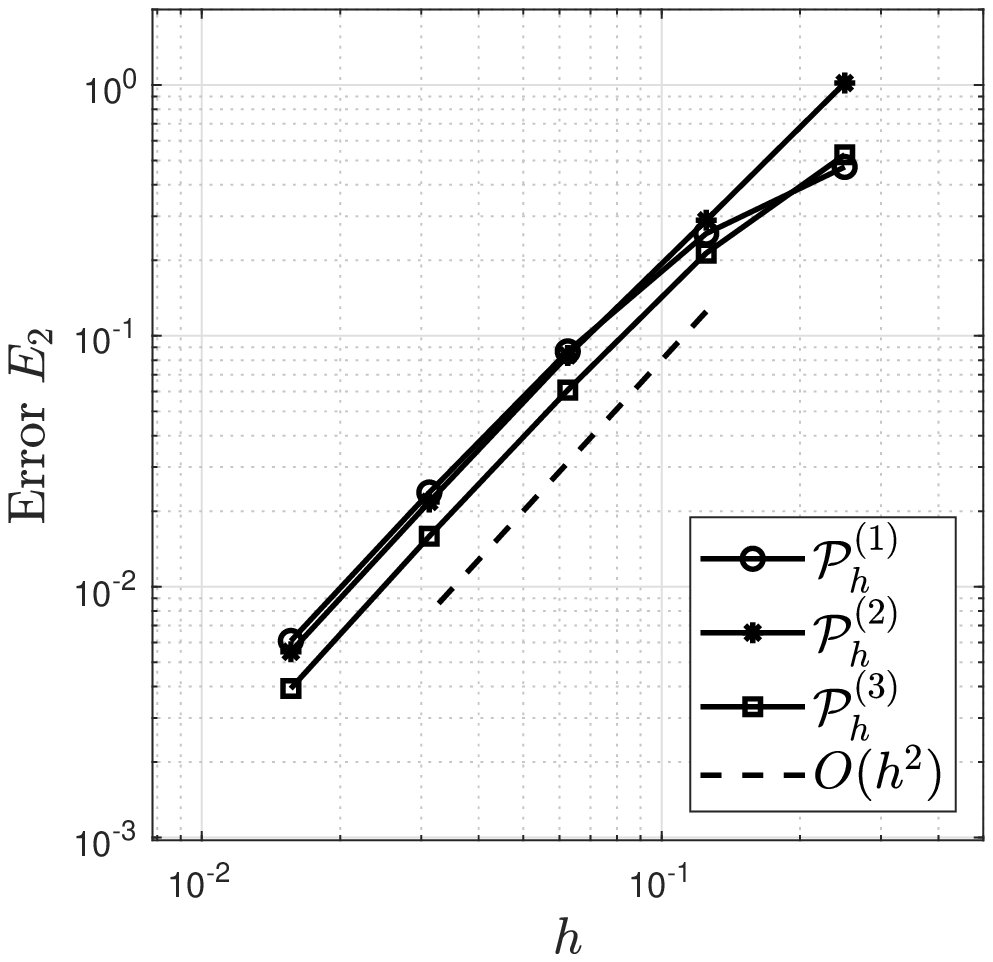}
\caption{The error curves $E_e$ (left) and $E_2$ (right) of test case 2 with $\lambda = 10^4$.}
\label{fig:errorKS2}
\end{figure}

We next implement the method \eqref{eqn:DiscreteProbKS} and compute the errors as above. In \Crefrange{fig:errorKS1}{fig:errorKS2}, we present the error curves versus $h$ for different values of $\lambda$. As shown in these figures, we see that the convergence order of the errors $E_e$ and $E_2$ are $O(h)$ and $O(h^2)$, respectively. Moreover, the convergence order is maintained in the nearly incompressible case ($\lambda = 10^4$). These results are consistent with the convergence rate predicted by the analysis in \Cref{thm:KSVEMconv}. 

\section{Conclusions}

We proposed two kinds of lowest-order virtual element methods for the linear elasticity problem. For the first one, we used the lowest-order virtual element method with a stabilizing term. This method can be seen as a modification of the Crouzeix-Raviart nonconforming finite element method as suggested in \cite{MR1972650} to the virtual element method. For the second one, we studied Kouhia-Stenberg type virtual element space, which consists of the conforming virtual element space for one component of the displacement vector and the nonconforming virtual element space for the other. This method can be seen as an extension of the Kouhia-Stenberg finite element method suggested in \cite{MR1343077} to the virtual element method. We proved that proposed methods have the optimal convergence of the numerical approximation to the dispacement vector field, and that the convergence is locking-free, that is, is stable with respect to $\lambda$. Finally, we present some numerical experiments that confirm the theoretical results. 

\bibliographystyle{siam}
\bibliography{refs}

\end{document}